\documentclass[aps,preprint,a4paper,nofootinbib,floatfix]{revtex4}
\usepackage{graphicx,color}
\usepackage{amsmath,amssymb,amsfonts,amsthm,amscd,bm}
\usepackage{booktabs}
\usepackage{listings}
\usepackage{float}

\def\beq{\begin{equation}}
\def\eeq{\end{equation}}

\def\tilde{\widetilde}

\def\[{\left[}
\def\]{\right]}
\def\({\left(}      
\def\){\right)}

\def\inv#1{\dfrac{1}{#1}}

\def\Li{{\rm Li}}

\def\tzero{t_\bullet}

\def\SL{SL_2\(\mathbb{Z}\)}
\def\CD{\mathcal{D}}
\def\sig{\nu}

\renewcommand{\arraystretch}{1}
\renewcommand{\arraycolsep}{5pt}
\makeatletter
\renewcommand*\env@matrix[1][\arraystretch]{%
  \edef\arraystretch{#1}%
  \hskip -\arraycolsep
  \let\@ifnextchar\new@ifnextchar
  \array{*\c@MaxMatrixCols c}}
\makeatother
\makeatletter
\g@addto@macro\bfseries{\boldmath}
\makeatother

\theoremstyle{plain}

\theoremstyle{remark}
\newtheorem{remark}{Remark}

\begin{document}

\title{
Transcendental equations satisfied by the  individual zeros of 
Riemann $\zeta$, Dirichlet and modular $L$-functions}

\author{Guilherme Fran\c ca\footnote{guifranca@gmail.com}
and  Andr\'e  LeClair\footnote{andre.leclair@gmail.com}}

\affiliation{Cornell University, Physics Department, Ithaca, NY 14850}


\begin{abstract}
We consider the non-trivial zeros of the Riemann $\zeta$-function
and two classes of $L$-functions; Dirichlet $L$-functions
and those based on level one modular forms.
We show that there are an infinite number  of zeros on the critical 
line in one-to-one correspondence with the zeros of the cosine function,
and thus  enumerated by an integer $n$. From this it follows that 
the ordinate  of the $n$-th  zero satisfies  a transcendental 
equation that depends only on $n$.
Under weak assumptions, we show that the number of solutions of this  
equation already saturates the counting formula on the entire critical 
strip.
We compute numerical solutions of these  transcendental equations and also 
its asymptotic limit of large ordinate. The starting point is an explicit 
formula, yielding an approximate solution for the ordinates of the zeros in 
terms of the Lambert $W$-function.
Our approach is a novel and simple method,
that takes into account $\arg L$, to numerically compute
non-trivial zeros of $L$-functions.  
The method is surprisingly accurate, fast and easy to implement.
Employing these numerical solutions, in particular
for the $\zeta$-function, we verify that the leading order asymptotic 
expansion is accurate enough to numerically support Montgomery's  
and Odlyzko's pair 
correlation conjectures, and also to reconstruct the prime number counting
function. Furthermore, the numerical solutions of the exact 
transcendental equation can determine the ordinates of the zeros
to any desired accuracy.
We also study in detail Dirichlet $L$-functions and  
the $L$-function for the modular form based on  
the Ramanujan $\tau$-function, which is closely related to the bosonic 
string partition function.
\end{abstract}

\maketitle
\tableofcontents

\section{Introduction}

Riemann's major contribution to number theory was an explicit formula 
for the arithmetic function $\pi (x)$, which counts the number of primes 
less than $x$, in terms of an infinite sum over the non-trivial zeros of 
the $\zeta (s) $ function,  i.e. roots  $\rho$ of the equation 
$\zeta (\rho) =0$ on the \emph{critical strip} 
$0\leq \Re (s)  \leq 1$ \cite{Riemann}.
It was later proven by Hadamard and de la Vall\'ee Poussin that there are 
no zeros on the line $\Re (s) =1$,  which in turn  proved the prime number 
theorem, $\pi (x)  \approx  \Li (x)$.  
(See section \ref{sec:prime} for a review.)
Hardy proved that there are an infinite number of zeros on the 
\emph{critical line} $\Re (s) = \tfrac{1}{2}$. 
The \emph{Riemann Hypothesis} (RH) was Riemann's statement in his
seminal eight-page paper \cite{Riemann} that all non-trivial zeros
have $\Re(\rho) = \tfrac{1}{2}$. In his own words, concerning the 
roots $t$ to the equation $\zeta(\tfrac{1}{2} + i t)=0$,
\begin{center}
\emph{``\ldots it is very likely that all roots are real. 
One would of course
like to have a rigorous proof of this, but I have put aside the search
for such a proof after some fleeting vain attempts \ldots''.}
\end{center}
Despite strong numerical evidence of its validity,  it remains unproven 
to this day. Many important mathematical results were proven assuming the RH, 
so it is a cornerstone of fundamental mathematics. Some excellent 
introductions to the RH are \cite{Edwards,Sarnak,Bombieri,Conrey}.

Riemann also gave an estimate $N(T)$, given by \eqref{riemann_counting} 
but without the $S(T)$ term, for the average number of zeros on the entire 
critical strip with $0 < \Im(\rho) < T$.
This formula was later proven by von Mangoldt, but it has never been proven 
to be valid on the critical line,  as explicitly stated in 
Edward's book \cite{Edwards}.
Denoting the number of zeros on the critical line up to height $T$ 
by $N_0(T)$, Hardy and
Littlewood proved that $N_0(T) > C \, T$. Selberg improved this result
stating that $N_0(T) > C \, T\log T$ for very small $C$. 
Levinson \cite{Levinson} demonstrated that 
$N_0(T) \ge C N(T)$ where $C = \tfrac{1}{3}$, which was
further improved by  Conrey \cite{Conrey2} who obtained $C = \tfrac{2}{5}$.
Further improvements on this last result are in \cite{Bui,Feng}.
Obviously, if the RH is true we must have
$N_0(T)=N(T)$. These statements
are described in \cite{Edwards,Titchmarsh}.

The RH is formulated as a problem in pure mathematics,  rather than physics,  
however it has interesting
connections with different areas of physics such as  quantum mechanics, 
quantum chaos,   and  in particular quantum statistical physics.
For an extensive review  we refer
the reader to \cite{Schumayer}. 
Although the work presented here does not intrinsically bring physics 
ideas  to bear 
on the problem,   
and is essentially pure mathematics,   it is worthwhile mentioning some ideas
on the RH that are based on physics,  even  if the purpose is  only to 
make contrasts with the present work.   
Julia \cite{Julia} and Spector \cite{Spector} proposed independently 
the free ``Riemann gas'' where the partition function is $\zeta$.
In \cite{Spector} supersymmetry and the Witten index were key ingredients.
A string
theory perspective on the RH is also possible.
Bakas and Bowick \cite{Bakas} considered an arithmetic gas to
construct a formula for
boson-parafermion equivalence using properties of $\zeta$. 
Examples of exactly solvable models were also discussed.
Spector \cite{Spector2} considered dualities in field theory that
are related to arithmetic functions. These are analogues of dualities
in string theory. He introduced the notion of partial
supersymmetry, leading to a formulation of parafermions of noninteger
order and found a bosonic analog of the Witten index.
These arithmetic quantum theories have a partition function related
to $\zeta$, and possess, like string theory, a Hagedorn temperature.
In \cite{Cacciatori} the RH is reformulated in terms of ultraviolet 
relations occurring in perturbative closed strings.
A connection between Gromov-Witten invariants, topological
string theory and Riemann zeros has also been motivated \cite{He}.
More ideas relating the RH to strings and geometry
can be found in \cite{Lapidus}.
A connection of the RH to quantum gases in low dimensions  was proposed in
\cite{AL}. 

The most prominent idea related to physics goes back to an old idea
of Hilbert-P\'olya.   Below,  we will describe and  study 
Montgomery's conjecture that the ordinates of non-trivial zeros of 
the $\zeta$-function satisfy the statistics of the 
\emph{Gaussian Unitary Ensemble} (GUE) \cite{Montgomery}.  The latter led
Berry to propose that the ordinates of the non-trivial zeros 
are eigenvalues of a chaotic hamiltonian \cite{Berry1}. 
Berry's work indicates interesting connections of the RH to
quantum chaos and was further explored in numerous papers. For instance, 
in  \cite{BerryKeating,BerryKeating2}  an analogy between 
the ordinates of the Riemann zeros and energy levels of a (unknown) quantum 
hermitian operator with chaotic dynamics was proposed.   
The classical counterpart of
such a hypothetical quantum system is associated to the 
hamiltonian $H = xp$. Based on this approach, a mapping 
between the Berry-Keating model and the Russian doll model of 
superconductivity was proposed \cite{Sierra}. 
This model is exactly solvable and has  a cyclic Renormalization Group.
In \cite{Sierra2} a generalization of the Berry-Keating model was 
considered  by adding an interaction term to the  hamiltonian.   
All these works focus on $N(T)$ and  carry out the a\-na\-ly\-sis on the 
critical line, i.e. they essentially assume the validity of the RH.   
Nevertheless,  
a number of interesting analytic results were obtained, emphasizing the 
important role of the fluctuating term in the counting formula $N(T)$,
namely the function $S(T)=\tfrac{1}{\pi}\arg\zeta\(\tfrac{1}{2}+iT\)$.
However, these works can only reproduce the smooth part of $N(T)$ 
through a semi-classical approach.
Instead of associating Riemann zeros to eigenstates of a 
quantum hamiltonian, as in the previously mentioned papers, the authors of
\cite{Bhaduri} focus on the scattering problem. They associate the 
smooth phase of the $\zeta$-function to the density of states of a quantum 
inverted harmonic oscillator.
In a related, but essentially different approach than Berry and Keating,
Connes used abstract mathematical objects called adeles. In this approach
there exists an operator playing the role of the hamiltonian,
which has a continuous spectrum, and the Riemann zeros correspond 
to missing spectral lines \cite{Connes}. Connes proposed a hypothetical 
trace formula which, if proved, can lead to a prove of the RH.
A dynamical system whose partition function is the $\zeta$-function
was also proposed \cite{Connes2}.
Unfortunately, thus far, a quantum mechanical hermitian operator 
whose spectrum yields the
non-trivial zeros has not yet been found. A quantum field theoretical
construction with a spectrum given by the Riemann zeros has also been pursued, 
although a free bosonic field theory with  a 
spectrum related to prime numbers is unlikely \cite{Svaiter,Andrade}, since
its path integral cannot be zeta-regularized. 
We will not be pursuing these ideas here,  rather,  the basis  of our work 
is a novel mathematical analysis of the original problem.     

$L$-functions are generalizations of the 
Riemann $\zeta$-function, the latter being the trivial case
\cite{Apostol}. In this paper we will consider two different classes of 
$L$-functions;  Dirichlet $L$-functions  and  $L$-functions 
associated with modular forms.
The former have applications primarily in multiplicative number theory,
whereas the latter in additive number theory.
These functions can be analytically continued to the entire 
(upper half) complex plane. The \emph{Generalized Riemann Hypothesis} (GRH)  
is the conjecture that all non-trivial zeros of Dirichlet $L$-functions 
and global $L$-functions in general lie on the critical line.
Much less is known about the zeros of $L$-functions in 
comparison with the $\zeta$-function,  however let us mention a few works.   
Selberg \cite{Selberg1} obtained the analog of Riemann-von Mangoldt 
counting formula \eqref{riemann_counting} for Dirichlet $L$-functions.
Based on this result, Fujii \cite{Fujii} gave an estimate for the number 
of zeros in  the critical strip with the ordinate between $[T, T+H]$. 
The distribution of low lying zeros of $L$-functions near and at the 
critical line was examined in \cite{Iwaniec}, assuming the GRH.
The statistics of the zeros,  i.e. the analog of the 
Montgomery-Odlyzko conjecture, were studied in \cite{Conrey3,Hughes,Bogomolny}.
It is also known that more than half of the non-trivial zeros 
of Dirichlet $L$-functions are on the critical line \cite{Conrey4}.
For a more detailed introduction to $L$-functions see \cite{Bombieri2}.

Besides the Dirichlet $L$-functions, there are more general constructions
of $L$-functions based on arithmetic and geometric objects, 
like varieties over number fields and modular forms \cite{Iwaniec2, Sarnak2}.
Some results for general $L$-functions are still conjectural. For instance,
it is not even clear if some $L$-functions can be analytically 
continued into a meromorphic function.   
We will only consider   the  additional $L$-functions
based on modular forms here.      
Thus the  $L$-functions considered in this paper have similar properties, 
namely, they possess  an Euler product,  
can be analytic continued into the (upper half) complex plane, 
except for possible poles at $z=0$ and $z=1$, 
and satisfy a non-trivial  functional 
equation.

Since it is well known that there are an  infinite number of zeros on 
the critical line for
the Riemann $\zeta$-function, if in some region of the critical strip 
one can show that the 
counting formula \eqref{riemann_counting} correctly counts the zeros on 
the critical line, then this proves the RH in this region of the strip.
It has been shown numerically that the first billion or so zeros all 
lie on the critical line \cite{deLune,Gourdon}, thus one can approach this
problem asymptotically. Such an analysis was carried out 
in \cite{RHLeclair},  where the main outcome  was an asymptotic 
transcendental equation for the ordinate of the $n$-th Riemann zero on the 
critical line.
The way in which this equation is derived shows that these zeros are in 
one-to-one correspondence with the zeros of the cosine function;
it is in this manner that the $n$-dependence arises.
In this paper we provide a more rigorous and through analysis of this result.
Moreover, we propose generalizations. We derive an \emph{exact} equation 
satisfied by the Riemann zeros on the critical line, where the above  mentioned 
asymptotic equation is obtained as a limit of large $n$. 
We also generalize these results to  Dirichlet
$L$-functions and to  $L$-functions related to modular forms. 
For all these classes of functions we obtain an exact equation for 
the ordinate of the $n$-th zero on the critical line. Since such an equation 
comes from a relation with the cosine function, its solutions 
can be automatically counted. We will argue that, under weak assumptions,
the number of solutions of the transcendental equation coincide with
the known counting formula for zeros on the entire  critical strip,  
i.e. $N_0 (T) = N(T)$. 

We organize our work as follows. In Section~\ref{sec:zeta_function} we 
derive an exact equation satisfied by each 
individual Riemann zero on the critical line. We discuss how
the number of its solutions can be the same as the counting formula 
on the entire critical strip.
In Section~\ref{sec:dirichlet} we follow the same analysis for
Dirichlet $L$-functions, and in Section~\ref{sec:modular} for $L$-functions
based on level one modular forms. In Section~\ref{sec:lambert} we derive a 
useful  approximation for the zeros expressed explicitly in terms of 
the Lambert $W$-function.
In Section~\ref{sec:davenport} we consider the counterexample
of Davenport-Heilbronn, which is known to violate the RH,  and discuss 
how the RH fails based on  the different  properties of our 
transcendental equation in 
comparison with previous cases.  
In Section~\ref{sec:numerical} we obtain numerical solutions to
the transcendental equation related to the Riemann $\zeta$-function.
We show that the leading order asymptotic 
approximation is accurate enough to reproduce the GUE statistics and the 
prime number counting function.
Furthermore, we show that solutions to the exact transcendental equation yield
highly accurate results,  up to $500$ digit accuracy or more if desired.   
In Section~\ref{sec:numerical_lfunc} we solve numerically the transcendental
equation related to Dirichlet $L$-functions, considering two explicit examples.
We also consider numerical solutions for $L$-functions based on modular 
forms, in particular for the $L$-function based on the 
Ramanujan $\tau$-function,  which is related to the bosonic string theory.  
Section~\ref{sec:conclusion} contains our concluding remarks.
In Appendix~\ref{sec:mathematica} we present the short 
Mathematica code we used to calculate the zeros for 
Dirichlet $L$-functions, some of which are
shown in Appendix~\ref{sec:tables}.

\section{Zeros of the Riemann $\zeta$-function}
\label{sec:zeta_function}

For simplicity we first consider the Riemann $\zeta$-function,  
which is the simplest
Dirichlet $L$-function. 
Moreover, we first consider the asymptotic equation 
\eqref{FinalTranscendence}, first proposed in \cite{RHLeclair}, since it 
involves more familiar functions.  However,
this asymptotic equation should  here be viewed as following straightforwardly 
from the new exact equation \eqref{exact_eq2}, presented later.

\subsection{Asymptotic equation satisfied by the $n$-th zero}
\label{sec:zeta_asymptotic}

Let us start with the completed Riemann zeta function defined by
\beq
\label{chidef}
\chi(s) \equiv \pi^{-s/2} \, \Gamma\(s/2\) \zeta(s)
\eeq
where $s = \sigma + i t$.
In quantum statistical physics, this function is the free energy of a 
gas of massless bosonic particles in $d$ spatial dimensions 
when $s=d+1$ \cite{AL},  up to the overall power of the 
temperature $T^{d+1}$.     
Under a ``modular'' transformation that exchanges one spatial 
coordinate with Euclidean time,  if one analytically  continues $d$,   
physical arguments shows that it must have 
the symmetry 
\beq
\label{chisym}
\chi\(s\) = \chi\(1-s\).
\eeq
This is the fundamental, and amazing, functional equation 
satisfied by the $\zeta$-function, which was proven by  Riemann using
only complex analysis.   
For several different ways of proving \eqref{chisym} see \cite{Titchmarsh}. 
Now consider Stirling's approximation
$\Gamma (s)  \simeq \sqrt{2 \pi}  s^{s - 1/2}e^{-s}$,
which is valid for large $t$. Under this
condition we also have
\beq
s^s = \exp\left\{ i\(t\log t + \dfrac{\pi \sigma}{2}\) + \sigma 
\log t - \dfrac{\pi t}{2} 
+ \sigma + O\(t^{-1}\) \right\}.
\eeq
Therefore, using the polar representation $\zeta = |\zeta| e^{i\arg\zeta}$ and
the above expansions,  we can write
$\chi(s) = A \, e^{i\theta}$ where
\begin{align}
A(\sigma,t) &= \sqrt{2\pi } \, \pi^{-\sigma/2} \(  
\dfrac{t}{2} \)^{(\sigma-1)/2} 
e^{- \pi t /4} |\zeta (\sigma + i t)|\(1+O\(s^{-1}\)\) ,  
\label{A_assymp} \\
\theta(\sigma, t) &= \dfrac{t}{2} \log \( \dfrac{t}{2 \pi e} \)  + 
\dfrac{\pi}{4}(\sigma-1) + \arg \zeta(\sigma + i t) + O\(t^{-1}\) .
\label{theta_assymp}
\end{align}
The above approximation is very accurate. For $t$ as low as $100$, 
it evaluates $\chi\(\tfrac{1}{2}  + i t\)$ correctly to one part in $10^6$.
Above, we are assuming $t > 0$. The results for $t<0$ follows trivially
from the relation $\overline{\chi(s)} = \chi(\overline{s})$.

Now let $\rho = \sigma+it$ be a Riemann zero. Then $\arg \zeta(\rho)$ can be 
defined by the limit
\beq \label{deltadef} 
\arg\zeta\(\rho\) \equiv  \lim_{\delta \to 0^+} \arg \zeta\(\sigma+\delta+it\).
\eeq
For reasons that are explained below,  
it is important  that ${0 < \delta \ll 1}$. 
This limit in general is not zero.
For instance, for the first Riemann zero given by
$\rho_1 \approx \tfrac{1}{2}+i\, 14.1347$, we have 
$\arg \zeta\(\rho_1\) \approx 0.157873919880941213041945$.
On the critical line $s=\tfrac{1}{2}+it$, if $t$ does not
correspond to the imaginary part of a zero, the well-known function 
$S(t) = \tfrac{1}{\pi}\arg\zeta\(\tfrac{1}{2}+it\)$
is defined by continuous variation along the straight 
lines starting from $2$, then up to $2+it$ and finally to
$\tfrac{1}{2}+it$, where $\arg\zeta(2)=0$. Assuming the RH, the current 
best bound is 
$|S(t)|\le\(\tfrac{1}{2}+o(1)\)\tfrac{\log t}{\log\log t}$ for 
$t\to \infty$, proven by Goldston and Gonek \cite{Goldston}. 
On a zero, the more  standard way to define this term is through the limit
$S(\rho) = \tfrac{1}{2} 
\lim_{\epsilon\to0} \( S\(\rho+i\epsilon\)+S\(\rho-i\epsilon\) \)$.
We have checked numerically that for several zeros on the line, our definition
\eqref{deltadef} gives the same answer as this standard approach,  
and also agrees
with the standard definition of $S(t)$ where $t$ is not the ordinate of 
a zero.   

From \eqref{chidef} we have $\overline{\chi(s)} = \chi\(\overline{s}\)$, 
which implies that
$A(\sigma,-t)=A(\sigma,t)$ and $\theta(\sigma,-t)=-\theta(\sigma,t)$. 
Denoting $\chi\(1-s\) = A' \, e^{-i\theta'}$ we then have
\beq\label{aatt}
A'(\sigma,t)=A(1-\sigma,t), \qquad \theta'(\sigma,t)=\theta(1-\sigma,t).
\eeq
From \eqref{chisym} we also have $|\chi(s)| = |\chi(1-s)|$, 
therefore $A(\sigma,t) = A'(\sigma,t)$ for \emph{any} $s$ on the 
critical strip.

Let us now approach  a zero
$\rho = \sigma + it$ through the $\delta \to 0^+$ limit. 
From \eqref{chidef} it follows that $\zeta(s)$ and $\chi(s)$ have
the same zeros on the critical strip, so it is enough
to consider the zeros of $\chi(s)$. From \eqref{chisym} we see 
that if $\rho$ is a zero so is $1-\rho$. Then we clearly have\footnote{
The linear combination in \eqref{sumchi}  was chosen 
to be manifestly symmetric under $s\to 1-s$. Had we taken a different 
linear combination in \eqref{sumchi}, 
then $B = e^{i\theta} + b \, e^{-i\theta'}$ for some constant $b$.
Setting the real and imaginary parts of $B$ to zero gives the two equations
$\cos \theta + b\cos \theta' =0$ and $\sin\theta - b  \sin \theta'=0$.    
Summing the squares of these equations one obtains  
$\cos (\theta + \theta')=-(b+1/b)/2$. However, since $b+1/b >1$,  there 
are no solutions except for $b=1$.}
\beq \label{sumchi}
\lim_{\delta\to 0^{+}}\[ \chi(\rho+\delta) + \chi(1-\rho-\delta) \] = 
\lim_{\delta \to 0^{+}} A(\sigma+ \delta,t) B(\sigma+\delta, t) = 0, 
\eeq
where we have defined
\beq\label{Bdef1}
B(\sigma,t) \equiv e^{i \theta(\sigma,t) }  +  e^{-i \theta'(\sigma,t)}.
\eeq
The second equality in \eqref{sumchi} follows from $A=A'$. 
Then, in the limit $\delta \to 0^{+}$, a zero corresponds 
to $A=0$, $B=0$ or both.   
They can simultaneously be zero since they are not independent.
If $B=0$ then $A=0$, since $A \propto |\zeta |$. However, the converse 
is not necessarily true.
 In order to be more rigorous,  one should 
consider the limits $\delta \to 0^+$ separately in $A$ verses $B$;  
below we will
consider taking the limit in $B$ first.   

The  non-trivial behavior of $A$ is mostly dictated by $|\zeta|$. On
the other hand there is much more structure in $B$ since it contains the
phases of $\chi(s)$ and $\chi(1-s)$. It describes oscillations
on the complex plane and involves $t \log t$ and $\zeta$ itself.
Thus let us consider $B=0$.    We will provide ample evidence that all zeros
are characterized by this equation.    
The general solution of $B=0$   is given by 
\beq\label{gen_sol}
\theta + \theta' = (2n+1)\pi,
\eeq
which are a family of curves $t(\sigma)$.
However, since $\chi(s)$ is analytic on the critical strip, 
we know that the zeros must be isolated  points rather than curves, 
thus this general solution must be restricted. 
Let us choose the particular solution
\beq \label{particular_sol}
\theta = \theta', \qquad \lim_{\delta\to 0^{+}}\cos\theta = 0.
\eeq
On the critical line $\sigma=\tfrac{1}{2}$, 
from \eqref{aatt} we have that the first equation in 
\eqref{particular_sol} is 
already satisfied. Then from the second equation in 
\eqref{particular_sol} we obtain
$\lim_{\delta\to0^+}\theta\(\tfrac{1}{2}+\delta,t\) = \(n+\tfrac{1}{2}\)\pi$ 
for $n=0,\pm 1, \pm 2, \dotsc$, hence
\beq
\label{almost_final_zeta}
n = \dfrac{t}{2\pi}\log\(\dfrac{t}{2\pi e}\) -\dfrac{5}{8}
+ \lim_{\delta\to0^{+}}\dfrac{1}{\pi}\arg \zeta\(\tfrac{1}{2}+\delta+i t\).
\eeq
A closer inspection shows that the right hand side of
\eqref{almost_final_zeta} has a minimum in the interval $(-2, -1)$, thus $n$ 
is bounded from below, i.e. $n \ge -1$. 
Establishing the \emph{convention} that zeros 
are labeled by positive integers, $\rho_n = \tfrac{1}{2}+i t_n$ where
$n=1,2,\dotsc$, we must replace $n \to n - 2$ in \eqref{almost_final_zeta}. 
Therefore, the imaginary
parts of these zeros satisfy the transcendental equation
\beq \label{FinalTranscendence} 
 \dfrac{t_n}{2 \pi}  \log \( \dfrac{t_n }{2 \pi e} \)   
+ \lim_{\delta \to 0^{+}}  \dfrac{1}{\pi} 
\arg  \zeta \( \tfrac{1}{2}+ \delta + i t_n  \) = n - \dfrac{11}{8}.
\eeq
In short, we have shown that, asymptotically, there are an  infinite 
number of  zeros on the critical line whose ordinates can be determined by 
solving \eqref{FinalTranscendence}. This equation determines
the zeros on the upper half of the critical line. The zeros on the lower
half are symmetrically distributed; if 
$\rho_n = \tfrac{1}{2}+it_n$ is a zero, so is 
$\overline{\rho}_n = \tfrac{1}{2}-it_n$.

The left hand side of \eqref{FinalTranscendence} is a monotonically  
increasing function of $t$, and the leading term is a smooth function. 
This is clear since the same terms appear in the staircase function 
$N(T)$, equation \eqref{riemann_counting}; see also Remark \ref{arg_term}.   
Possible discontinuities can only come from 
$\tfrac{1}{\pi}\arg\zeta\(\tfrac{1}{2}+it\)$, and in fact, it has a jump 
discontinuity  whenever $t$ corresponds to 
the ordinate of a zero on or off the critical line. 
However, if
$\lim_{\delta \to 0^+}\arg\zeta\(\tfrac{1}{2}+\delta+it\)$ is 
well-defined for every $t$, then the left hand side of equation 
\eqref{FinalTranscendence} is well-defined for any $t$,  and due to 
its monotonicity, there  must be a 
unique solution for every $n$.
Under this assumption, the number of solutions of equation 
\eqref{FinalTranscendence}, up to height $T$, is given by
\beq \label{counting2}
N_0(T) = \dfrac{T}{2 \pi} \log \( \frac{T}{2 \pi e} \) + \frac{7}{8}  + 
\inv{\pi} \arg \zeta \( \tfrac{1}{2} + i T \) + O\(T^{-1}\).
\eeq
This is so because the zeros are already numbered in 
\eqref{FinalTranscendence}.  
Thus we can replace $n \to N_0 + \tfrac{1}{2}$ and $t_n \to T$, such that the 
jumps correspond to integer values. In this way $T$ will not correspond 
to the ordinate of a zero and $\delta$ can be eliminated.
In summary,   $N_0 (T)$ in \eqref{counting2}  counts the solutions to the
equation \eqref{FinalTranscendence} for zeros on the critical line,   
assuming there is a solution for every $n$,
and without assuming the RH.   

Using Cauchy's argument principle it is known that one can derive 
the Riemann-von Mangoldt 
formula, which gives the number of zeros inside the critical
strip with $0< \Im(\rho) < T$.
This formula is given by \cite{Edwards, Titchmarsh}
\beq \label{riemann_counting}
N(T) = \dfrac{T}{2 \pi} \log \( \frac{T}{2 \pi e} \) + \frac{7}{8}  + 
S(T) + O\(T^{-1}\)
\eeq
where $S(T) = \tfrac{1}{\pi}\arg\zeta\(\tfrac{1}{2}+iT\)$.
The above formula without the $S(T)$ term was already
in Riemann's paper \cite{Riemann}.   
Note that it has the same form as the counting formula  on 
the critical line that we have just found, equation \eqref{counting2}.   
Thus, under the assumptions we have described,  we conclude that 
$N_0(T) \simeq N(T)$,  at least asymptotically for now.   In the next
section we will present the exact version.   
This means that our particular solution \eqref{particular_sol}, leading
to equation \eqref{FinalTranscendence}, already saturates the counting 
formula on the entire critical strip and there are no additional 
zeros from $A=0$ in \eqref{sumchi}, nor from the general solution 
\eqref{gen_sol}. This strongly suggests that
\eqref{FinalTranscendence} describes all non-trivial zeros of $\zeta$, 
which must then lie on the critical line.    
We emphasize that we have {\it not}   
assumed the RH in the above arguments.    
In Section~\ref{sec:conclusion} 
we will  summarize the assumptions  which lead  to the  exact version of the 
equation \eqref{FinalTranscendence} described below  and
reiterate  its implications for the RH.

\subsection{Exact equation for the $n$-th zero}
\label{sec:zeta_exact}

It is straightforward to repeat the above analysis without
considering an asymptotic expansion. The exact versions of
\eqref{A_assymp} and \eqref{theta_assymp} are
\begin{align}
\label{A_exact}
A(\sigma,t) &= \pi^{-\sigma/2} |\Gamma\(\tfrac{1}{2}(\sigma+it)\)| 
|\zeta(\sigma+it)|, \\
\label{theta_exact}
\theta(\sigma,t) &= 
\arg \Gamma\(\tfrac{1}{2}(\sigma+it)\) -\dfrac{t}{2}\log\pi + 
\arg \zeta(\sigma+it),
\end{align}
where again $\chi(s)=Ae^{i\theta}$ and $\chi(1-s)=A'e^{-i\theta'}$, with
$A'(\sigma,t)=A(1-\sigma,t)$ and $\theta'(\sigma,t) = \theta(1-\sigma,t)$. 
The zeros on the critical line $\sigma=\tfrac{1}{2}$ correspond to the  
particular solution $\theta=\theta'$ and 
$\lim_{\delta\to 0^+}\cos \theta =0$. Therefore 
we have 
$\lim_{\delta\to0^+}\theta\(\tfrac{1}{2}+\delta,t\) = \(n+\tfrac{1}{2}\)\pi$. 
Replacing $n \to n-2$, the imaginary parts of these zeros must 
satisfy the \emph{exact} equation
\beq\label{exact_eq}
\arg\Gamma\(\tfrac{1}{4}+\tfrac{i}{2}t_n\) - t_n \log\sqrt{\pi} 
+ \lim_{\delta \to 0^+} \arg\zeta\(\tfrac{1}{2}+\delta+it_n\) = 
\(n-\tfrac{3}{2}\)\pi.
\eeq
The Riemann-Siegel $\vartheta$ function is defined by
\beq
\label{riemann_siegel}
\vartheta(t)\equiv \arg\Gamma\(\tfrac{1}{4}+\tfrac{i}{2}t\) - t\log \sqrt{\pi},
\eeq
where $\arg\Gamma$ is defined such that this function is continuous and
$\vartheta(0)=0$.   Therefore, we conclude that there are an infinite  
number of zeros in the form
$\rho_n=\tfrac{1}{2}+it_n$, where $n=1,2,\dotsc$, whose imaginary
parts \emph{exactly} satisfy the following equation:
\beq\label{exact_eq2}
\vartheta(t_n) + 
\lim_{\delta\to 0^{+}}\arg\zeta\(\tfrac{1}{2}+\delta +it_n\) = 
\(n-\tfrac{3}{2}\)\pi.
\eeq
Expanding the $\Gamma$-function in \eqref{riemann_siegel} through 
Stirling's formula one obtains 
$\vartheta(t_n) = \tfrac{t_n}{2}\log\(\tfrac{t_n}{2\pi e}\)-\tfrac{\pi}{8}+
O\(1/t_n\)$, and recovers the asymptotic equation 
\eqref{FinalTranscendence} from \eqref{exact_eq2}.
Let us mention at this point  that our approach of considering zeros of $B=0$,  namely  
\eqref{gen_sol},   is
also able to reproduce the trivial zeros on the negative real
line,  and also zeros off of the critical line in the counterexample 
of Section~\ref{sec:davenport} \cite{Lectures}.     

Again, as discussed in the paragraph above equation 
\eqref{counting2}, the first
term in \eqref{exact_eq2} is smooth and the whole left hand side is
a monotonic increasing function. If 
$\lim_{\delta\to 0^{+}}\zeta\(\tfrac{1}{2}+\delta+it\)$ is well-defined
for every $t$, then equation \eqref{exact_eq2} must have a unique solution
for every $n$;  see also  
Section~\ref{sec:davenport}.  Under this assumption it is valid to replace 
$t_n \to T$ and $n \to N_0 + \tfrac{1}{2}$, so 
the number of solutions of \eqref{exact_eq2} is given by
\beq\label{counting2_exact}
N_0(T) = \dfrac{1}{\pi}\vartheta(T) + 1 +
\dfrac{1}{\pi}\arg\zeta\(\tfrac{1}{2}+iT\).
\eeq

The exact Backlund counting formula, which gives the number of zeros on
the entire critical strip with $0<\Im(\rho)<T$, is given by  the well-known 
formula \cite{Edwards}
\beq\label{backlund}
N(T) = \dfrac{1}{\pi}\vartheta(T) +  1 + S\(T\).
\eeq
Therefore, comparing with the exact counting formula on the entire 
critical strip \eqref{backlund}, we have $N_0(T) = N(T)$ exactly. 
This indicates that our par\-ti\-cu\-lar solution, 
leading to equation \eqref{exact_eq2}, captures all the zeros on the 
critical strip,
and they should all be on the critical line.

In summary, without assuming the RH, but \emph{under the assumption}  
that $\lim_{\delta\to0^+}\arg\zeta\(\tfrac{1}{2}+\delta+it\)$ exists
for every $t$, then \eqref{exact_eq2}  has a unique solution for 
every $n$.   If one ignores the $\arg \zeta$ term in the equation,
then the unique solution is expressed in terms of the
Lambert $W$-function;
see section \ref{sec:lambert}.   If there is indeed a unique solution for every $n$,  then
this leads to a $N_0 (T)$ which  saturates the counting formula for the entire 
critical strip,  and  this would establish the validity of the RH.  
Furthermore, it implies that
all non-trivial zeros are simple,  
as explained in Remark~\ref{simple_zeros}.    
Further related and clarifying  remarks,   based on a counterexample,   
are in Section~\ref{sec:davenport}.

\subsection{Further remarks}
\label{sec:zeta_remarks}

\begin{figure}
\centering
\begin{minipage}{.5\textwidth}
    \centering
    \includegraphics[width=1\linewidth]{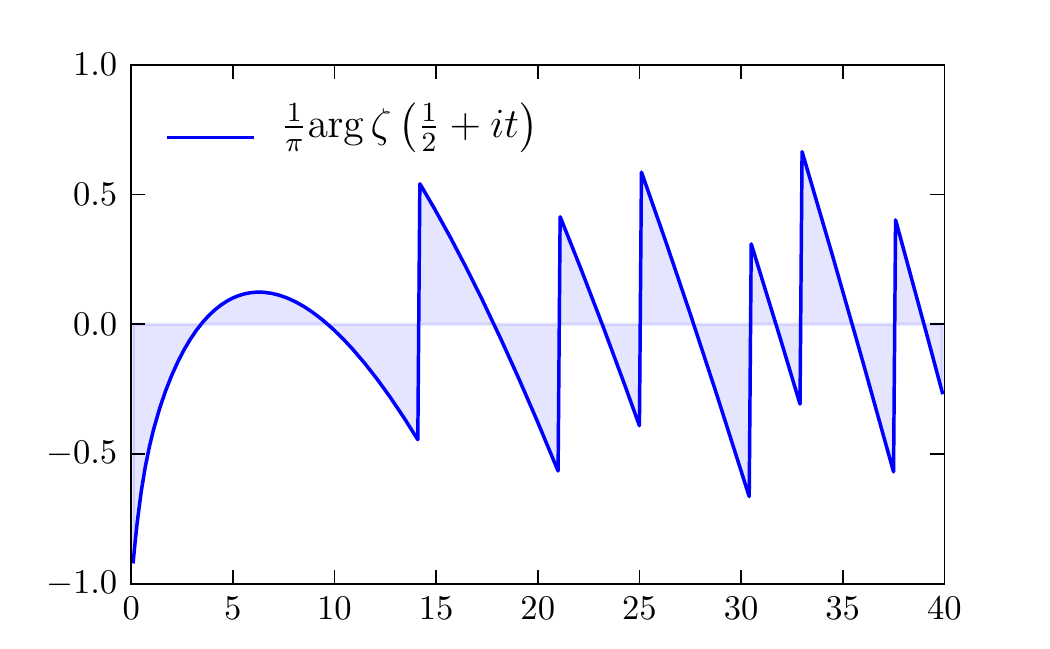}\\[-1em]
    \small{(a)}
\end{minipage}%
\begin{minipage}{.5\textwidth}
    \centering
    \includegraphics[width=1\linewidth]{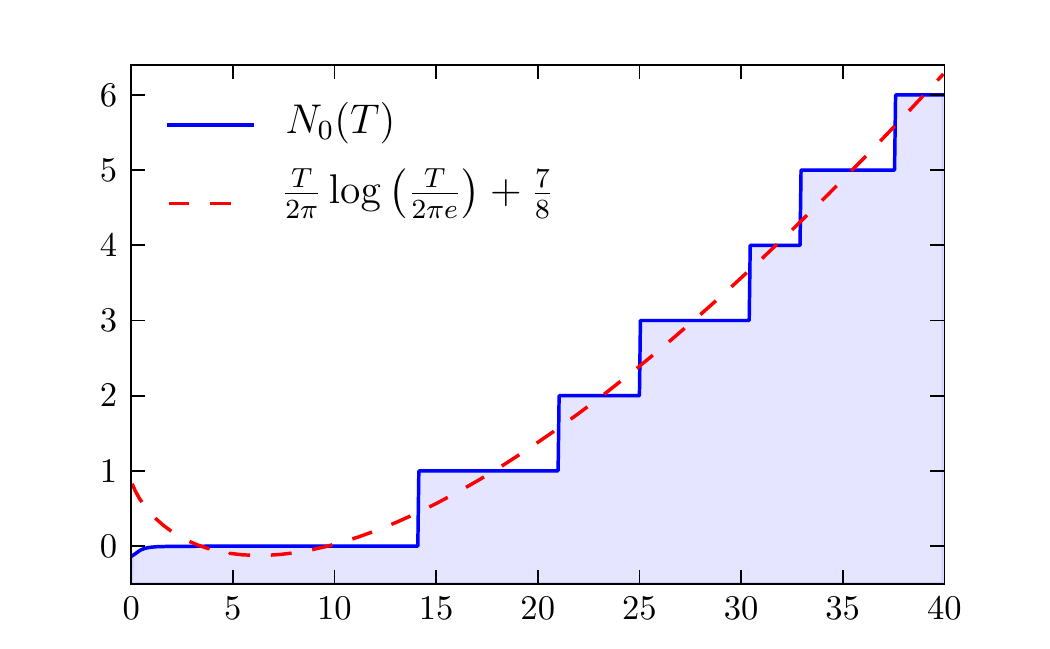}\\[-1em]
    \small{(b)}
\end{minipage}
\caption{(a) A plot of $\tfrac{1}{\pi}\arg \zeta\(\tfrac{1}{2} + i t\)$  
as a function of $t$ showing 
its rapid oscillation. The jumps occur on a Riemann zero. 
(b) The function $N_0(T)$ in \eqref{counting2}, which is indistinguishable 
from a manual counting of zeros.  The dashed line is the smooth part without
the $\arg \zeta$ term. }
\label{fig:arg_counting}
\end{figure}

\begin{remark}
\label{arg_term}  
The small shift by $\delta$ in equations \eqref{exact_eq2} or 
\eqref{FinalTranscendence} is essential
since it smooths out $S(t)  =  \tfrac{1}{\pi}\arg \zeta\(\tfrac{1}{2}+ i t\)$,
which is known to jump discontinuously at each Riemann zero.
As is well known, $S(t)$ is a piecewise continuous function
that rapidly 
oscillates around its average value, which is zero,  with 
discontinuous jumps, as shown in Figure~\ref{fig:arg_counting}a.
However, when $S(t)$ is added to the smooth part 
of $N(T)$ one obtains an accurate staircase function, which 
jumps by the multiplicity of the zero at 
the ordinate of each Riemann  zero; 
see Figure~\ref{fig:arg_counting}b.  
Note that $N(T)$ is necessarily a monotonically increasing function
since it is a counting formula.

One reason $\delta$ needs to be positive
in \eqref{exact_eq2} can be seen as follows. 
Near a simple zero $\rho_n$ we have
$\zeta (s)  \approx  \(s-\rho_n\) \zeta' \(\rho_n\)
= \(\delta +i \(t-t_n\)\) \zeta'\(\rho_n\)$.
This gives 
$\arg \zeta (s) \sim \arctan\((t-t_n)/\delta\)+ c$, where $c$ is a constant.
With $\delta > 0$  as one passes through a zero from below,  
$S(t)$ \emph{increases} by one, as it should based on its role in the counting
formula $N(T)$. On the other hand, if $\delta < 0$ then $S(t)$ would 
\emph{decrease} by one instead, which cannot be the 
case\footnote{Note added:  there 
are deeper reasons why $\delta$ has to be positive,  
described in our subsequent work
\cite{EulerProd},  which is discussed in the concluding section 
of the present article.}.

\end{remark}

\begin{remark}
\label{simple_zeros}
An important consequence of equation \eqref{exact_eq2} 
is that, again,
if it has a unique solution for every $n$,  then all non-trivial 
zeros are simple.
This essentially  follows from the fact that they 
are in one-to-one correspondence with the zeros of the cosine 
function \eqref{particular_sol}, which are simple.  To see this, 
let us suppose there is a \emph{double zero} at ordinate $\tzero$, i.e.
$t_{n+1} = t_n = \tzero$. Then subtracting the 
equation \eqref{exact_eq2} with $t_{n}$ from the corresponding equation
with $t_{n+1}$, we obtain
a contradiction, namely $0 = \pi$. Therefore, for $n\neq m$ the equation
\eqref{exact_eq2} implies $t_n \neq t_m$.

Now if we  actually assume that the zeros on the critical line are 
simple (which we have not),  
there is an easier  way to see 
that the zeros correspond to $\cos \theta =0$. On the critical
line $s = \tfrac{1}{2} + i t$, the functional equation 
\eqref{chisym} implies $\chi (s)$ is real, thus for $t$ not the ordinate
of a zero, $\sin\theta = 0$ and $\cos\theta = \pm 1$.
Thus $\cos \theta $ is a discontinuous function. 
Now let $\tzero$ be the ordinate of a simple zero. Then close
to such a zero we define 
\beq
c(t) \equiv  \frac{ \chi (\tfrac{1}{2} + i t ) }{|\chi (\tfrac{1}{2} + i t)|}  
\sim  \frac{ t-\tzero}{|t - \tzero|}.
\eeq
For $t > \tzero$ then $c(t)=1$, 
and for $t<\tzero$ then  $c(t) = -1$.  
Thus $c(t)$ is discontinuous precisely at the zero.
In the above polar representation, formally 
$c(t) =  \cos \theta \(\tfrac{1}{2}, t\)$. 
Therefore, by identifying zeros as the solutions to $\cos \theta =0$,
we are simply defining  the value of the function  $c(t)$ at 
the discontinuity as  $c(\tzero)=0$.
\end{remark}

\begin{remark}
\label{newzeta}
Let us introduce another function 
$\zeta(s) \to \widetilde{\zeta}(s) = f(s) \zeta(s)$ that also 
satisfies the functional equation \eqref{chisym}, 
i.e. $\tilde{\chi}(s)=\tilde{\chi}(1-s)$, but has zeros 
off of the critical line due to the zeros of $f(s)$. In such a case 
the corresponding functional equation will hold if and only 
if $f(s) = f(1-s)$ for any $s$, and this is a trivial condition 
on $f(s)$ which could have been canceled in the first place. 
Moreover, if $f(s)$ and $\zeta(s)$ have different zeros, the analog of 
equation \eqref{sumchi} has a  factor
$f(s)$, i.e. $\tilde{\chi}(\rho+\delta) + \tilde{\chi}(1-\rho-\delta) = 
f(\rho+\delta)\[\chi(\rho+\delta)+\chi(1-\rho-\delta)\]=0$,  
implying \eqref{sumchi} again  where $\chi(s)$ is the original \eqref{chidef}. 
Therefore, the previous analysis eliminates $f(s)$
automatically and only finds the zeros of $\chi(s)$.
The analysis is  non-trivial precisely  because $\zeta(s)$ satisfies 
the functional equation but $\zeta(s) \ne \zeta(1-s)$. Furthermore, 
it is a well-known theorem that the only function which satisfies the 
functional equation \eqref{chisym} and has the same characteristics 
of $\zeta(s)$, is $\zeta(s)$ itself. 
In other words, if $\widetilde{\zeta}(s)$ is required to have  the 
same properties of $\zeta(s)$, then $\widetilde{\zeta}(s) = C \, \zeta(s)$ 
where $C$ is a constant \cite[pg. 31]{Titchmarsh}.
\end{remark}

\begin{remark}
\label{counting}
Although equations \eqref{exact_eq2} and \eqref{backlund} have an obvious 
resemblance, it is impossible to derive the former from the later, 
since the later is just a counting formula valid on the entire strip, 
and it is assumed that $T$ is \emph{not} the ordinate of a zero.
Moreover,  such a derivation would require the assumption of the validity 
of the RH and the simplicity of the zeros,
contrary to our approach, where we derived equations \eqref{exact_eq2}
and \eqref{FinalTranscendence} directly on the critical line, without assuming
the RH,  nor the known counting formula $N(T)$.   
Despite our best efforts, we were not able to find 
equations \eqref{FinalTranscendence} and \eqref{exact_eq2} 
in the literature. 
Furthermore,  the counting formulas \eqref{counting2} and \eqref{backlund}
have  never been proven  to be valid on the critical line \cite{Edwards}.
\end{remark}

\begin{remark} 
\label{circular}
One may object that our basic equation 
\eqref{exact_eq2}  involves  $\zeta(s)$ itself and this is 
somehow circular.
This is not a valid counter-argument.
First of all, $\arg \zeta$ already appears in the counting function $N(T)$.  
Secondly, the equation \eqref{exact_eq2} is a much more 
detailed equation than simply $\zeta (s )= 0$,  which has an infinite 
number of solutions,  in contrast with  \eqref{exact_eq2}  which 
for each $n$,  as we have argued,   
has  a unique solution corresponding to the $n$-th zero.
Also, there are well known ways to calculate $\arg \zeta$,  
for example from an integral representation or
a convergent series \cite{Borwein}. 
\end{remark}

\section{Zeros of Dirichlet $L$-functions}
\label{sec:dirichlet}

\subsection{Some properties of Dirichlet $L$-functions }  

We now consider the generalization of the previous results to Dirichlet
$L$-functions.  Let us  first  introduce the basic ingredients and definitions  
regarding this class of functions, which are all well known \cite{Apostol}.
Dirichlet $L$-series are defined as 
\beq
\label{Ldef}
L(s, \chi) = \sum_{n=1}^\infty \frac{\chi (n)}{n^s}
\eeq
for $\Re(s) > 1$, 
where the arithmetic function $\chi(n)$ is a Dirichlet character.  
They  can all be analytically continued to the entire 
complex plane, except possibly for a simple pole at $s=1$ when
$\chi$ is principal, and are then referred to as  Dirichlet $L$-functions. 

There are an infinite number of distinct Dirichlet 
characters which are primarily characterized by their modulus $k$,
which determines their periodicity. They can
be defined axiomatically,  which leads to specific properties,
some of which we now describe.
Consider a Dirichlet character $\chi$ mod $k$,  and let 
the symbol $(n,k)$ denote the greatest  common divisor of the
two integers $n$ and $k$. Then $\chi$  has the 
following properties:
\begin{enumerate}
\item $\chi(n+k) = \chi (n)$.
\item $\chi(1) = 1$ and  $\chi(0) = 0$.
\item $\chi(n \, m) = \chi (n) \chi (m)$.
\item $\chi(n) = 0$ if  $(n, k) > 1$ 
and $\chi(n) \ne 0$ if $(n, k) = 1$.
\item \label{root1} If $(n, k) = 1$  then  $\chi(n)^{\varphi(k)} = 1$, where
$\varphi(k)$  is the Euler totient arithmetic function. This implies  
that $\chi(n)$ are roots of unity.
\item If $\chi$ is a Dirichlet  character so is the 
complex conjugate $\overline{\chi}$.
\end{enumerate}
For a given modulus $k$ there are $\varphi(k)$ distinct Dirichlet 
characters,  which essentially follows from Property~\ref{root1} above.
They can thus be labeled as $\chi_{k,j}$ where $j= 1, 2,\dotsc, \varphi(k)$ 
denotes an arbitrary ordering. 
If $k=1$ we have the \emph{trivial} character where
$\chi(n)=1$ for every $n$, and \eqref{Ldef} reduces to the Riemann 
$\zeta$-function.
The \emph{principal } character,  usually denoted by $\chi_1$,   
is defined as 
$\chi_1(n) = 1$ if $(n,k) = 1$ and zero otherwise. In the above notation 
the principal character is always $\chi_{k,1}$.

Characters can be classified as \emph{primitive} or \emph{non-primitive}.
Consider the Gauss sum
\beq\label{tau}
G(\chi) = \sum_{m=1}^{k}\chi(m)e^{2\pi i m / k}.
\eeq
If the character $\chi$ mod $k$ is primitive, 
then $|G(\chi)|^2 = k$. This is no longer valid for a non-primitive character.
Consider a non-primitive character $\tilde{\chi}$ 
mod $\tilde{k}$. Then it can be expressed in terms of
a primitive character of smaller modulus
as $\tilde{\chi}(n) = \tilde{\chi_1}(n) \chi(n)$, where $\tilde{\chi_1}$ is
the principal character mod $\tilde{k}$ and $\chi$ is a primitive
character mod $k < \tilde{k}$, where $k$ is a divisor of $\tilde{k}$. More
precisely, $k$ must be the \emph{conductor} of $\tilde{\chi}$ 
(see \cite{Apostol} for further details).
In this case the two $L$-functions are related as
$L(s, \tilde{\chi}) = L(s,\chi) \Pi_{p | \tilde{k}}\(1 - \chi(p)/p^s\)$. 
Thus $L(s,\tilde{\chi})$ has the same zeros as $L(s,\chi)$.
Therefore, it suffices to consider 
primitive characters, and we will henceforth do so.  

We will need the functional equation satisfied by $L(s,\chi)$.   
Let $\chi$ be a \emph{primitive} character. 
Define its \emph{order} $a$ such that 
\beq\label{order}
a \equiv \begin{cases}
1 &\mbox{if $\chi(-1)=-1$ (odd)} \\
0 &\mbox{if $\chi(-1)=1$ (even)}
\end{cases}.
\eeq
Let us define the  entire function 
\beq
\label{Lambda}
\Lambda(s,\chi) \equiv \( \dfrac{k}{\pi} \)^{\tfrac{s+a}{2}}   
\Gamma\(\dfrac{s+a}{2}\) L(s, \chi).   
\eeq
Then $\Lambda$  satisfies the following 
well-known functional equation,
only valid for primitive characters \cite{Apostol}:
\beq
\label{FELambda}
\Lambda (s, \chi) =  \dfrac{ i^{-a}  \, G(\chi)}{\sqrt{k}} \,
\Lambda (1-s, \overline{\chi}).
\eeq

\subsection{Exact equation for the $n$-th zero}
\label{sec:dirichlet_equation}

For a primitive character, since $|G(\chi)| = \sqrt{k}$,  
the factor on the right hand side of \eqref{FELambda} is a phase. 
It is  thus possible to obtain a more symmetric form  
through a new function defined as
\beq
\label{xi}
\xi(s, \chi) \equiv \dfrac{i^{a/2} \, k^{1/4} }{\sqrt{G\(\chi\)}} \, \, 
\Lambda (s, \chi). 
\eeq
It  then satisfies 
\beq
\label{FExi}
\xi(s, \chi)  = \overline{\xi}(1-s , \chi ) \equiv   
\overline{\xi(1-\overline{s}, \chi)}.
\eeq
Above, the function $\overline{\xi}$  of $s$ is defined as the 
complex conjugation of
all coefficients that  define $\xi$,  
namely  $\chi$ and the $i^{a/2}$ factor, evaluated at a non-conjugated $s$.
   
Note that $\overline{\Lambda(s, \chi)} = 
\Lambda(\overline{s},\overline{\chi})$. Using the 
known result $G\(\overline{\chi}\) = \chi(-1) \overline{G(\chi)}$ we 
then conclude that
\beq\label{xiconju}
\overline{\xi(s,\chi)} = \xi\(\overline{s},\overline{\chi}\).
\eeq
This implies that if the character is real, when $\rho$ is a 
zero of $\xi$  so is $\overline{\rho}$, and one needs only to
consider $\rho$ with  
positive imaginary part. On the other hand if $\chi \neq \overline{\chi}$, 
then the zeros with negative imaginary part are different than 
$\overline{\rho}$.
For the trivial character where $k=1$ and $a=0$, implying 
$\chi(n) = 1$ for any $n$,  then $L(s,\chi)$  reduces to the
Riemann $\zeta(s)$ and  \eqref{FExi} yields the 
well-known functional equation \eqref{chisym}.

Let  $s=\sigma+i t$. Then 
the function \eqref{xi} can be written as
$\xi(s,\chi) = A e^{i \theta}$ where
\begin{align}
A(\sigma,t,\chi) & =
\(\dfrac{k}{\pi}\)^{\tfrac{\sigma+a}{2}}
\left|\Gamma\(\dfrac{\sigma+a+it}{2}\)\right| 
\left|L(\sigma+it,\chi)\right|, \label{Achi} \\
\theta(\sigma,t,\chi) &= \arg \Gamma\(\dfrac{\sigma+a+it}{2}\) - 
\dfrac{t}{2}\log \(\dfrac{\pi}{k}\) 
- \dfrac{1}{2}\arg G(\chi) \label{thetachi} \\
\nonumber & \qquad + \arg L(\sigma+it,\chi) + \dfrac{\pi a}{4}. 
\end{align}
From \eqref{xiconju} we have that
$A(\sigma,t,\chi) = A(\sigma,-t,\overline{\chi})$ 
and $\theta(\sigma,t,\chi) = -\theta(\sigma,-t,\overline{\chi})$.
Denoting $\overline{\xi}(1-s, \chi) = A' e^{-i \theta'}$
we then have 
$A'(\sigma,t,\chi) = A(1-\sigma,t, \chi)$ and 
$\theta'(\sigma,t, \chi) = \theta(1-\sigma,t,\chi)$.
Taking the modulus of \eqref{FExi} we also have that
$A(\sigma,t,\chi)=A'(\sigma,t,\chi)$ for any $s$.

On the critical strip, the functions $L(s,\chi)$ and $\xi(s,\chi)$ 
have the same zeros.
Thus on a zero
we clearly have
\beq\label{onzero}
\lim_{\delta\to0^+} 
\left\{ \xi (\rho+\delta, \chi)  +  
\overline{\xi} (1-\rho-\delta,\chi ) \right\} =0.
\eeq
Let us define
\beq
\label{Bdef}
B(\sigma,t,\chi) \equiv e^{i\theta (\sigma,t,\chi)} + 
e^{- i\theta' (\sigma,t,\chi)} .
\eeq
Since $A=A'$ everywhere,  from \eqref{onzero} we 
conclude that on a zero we have
\beq\label{theta_zero}
\lim_{\delta \to 0^+}  A(\sigma+ \delta ,t, \chi) 
B(\sigma+\delta ,t, \chi) = 0.
\eeq
As before, let us consider the particular solution
of $\lim_{\delta\to 0^+}B=0$ given by
\beq\label{particular}
\theta = \theta' , \qquad \lim_{\delta\to 0^{+}}\cos\theta = 0. 
\eeq
Let us define the function
\beq\label{RSgen}
\begin{split}
\vartheta_{k,a} (t) &\equiv  
\arg \Gamma \( \dfrac{1}{4} + \dfrac{a}{2} + i \, \dfrac{t}{2}  \)  
- \dfrac{t}{2}  \log \( \dfrac{\pi}{k} \) \\
& = \Im\[\log\Gamma\( \dfrac{1}{4} + \dfrac{a}{2} +i \,  \dfrac{t}{2} \)\]
- \dfrac{t}{2} \log \( \dfrac{\pi}{k} \).
\end{split}
\eeq
When $k=1$ and $a=0$, the function \eqref{RSgen} is just the usual 
Riemann-Siegel $\vartheta$ function \eqref{riemann_siegel}. 
Since the function $\log\Gamma$ has a complicated branch cut, one can
use the following series representation in \eqref{RSgen} \cite{Boros}
\beq
\log\Gamma(s) = -\gamma s - \log s -
\sum_{n=1}^{\infty}\left\{ \log\(1+\dfrac{s}{n}\) -\dfrac{s}{n}\right\},
\eeq
where $\gamma$ is the Euler-Mascheroni constant. Nevertheless,
most numerical packages already have the $\log\Gamma$ function implemented.

On the critical line $\sigma=\tfrac{1}{2}$ the first equation 
in \eqref{particular} is already satisfied. From the second equation we have 
$\lim_{\delta\to0^+}\theta\(\tfrac{1}{2}+\delta,t\)=\(n+\tfrac{1}{2}\)\pi$,  
therefore
\beq\label{almost_final}
\vartheta_{k,a}(t_n) + 
\lim_{\delta\to 0^{+}}\arg L\(\tfrac{1}{2}+\delta+it_n,\chi\) 
- \dfrac{\arg G}{2} + \dfrac{\pi a }{4}  
= \(n+\dfrac{1}{2} \)\pi .
\eeq
Analyzing the left hand side of \eqref{almost_final} we can see that it
has a minimum, thus we shift $n \to n - (n_0 + 1)$ for a given $n_0$, 
to label the zeros according to the convention that the first positive 
zero is labelled by $n=1$. 
Thus the upper half of the critical line will have the zeros
labelled by $n=1,2,\dotsc$ corresponding to positive $t_n$, while the 
lower half will have the negative values $t_n$ 
labelled by $n=0,-1,\dotsc$. The integer $n_0$ depends 
on $k$, $a$ and $\chi$, and should be chosen according to each specific case.
In the cases we analyze below $n_0=0$, whereas for the trivial 
character $n_0 =1$. 
Henceforth we will \emph{omit} the integer $n_0$ in the equations, 
since all cases analyzed in the following have $n_0=0$. 
Nevertheless, the reader
should bear in mind that for other cases, it may  be necessary
to replace $n \to n - n_0$ in the following equations.

In summary,  there are an  infinite number of  zeros on 
the critical line, i.e. in the  form $\rho_n = \tfrac{1}{2}+i t_n$, 
where for a given  
$n \in \mathbb{Z}$, the imaginary part $t_n$ is the solution of 
the equation
\beq\label{exact}
\vartheta_{k,a}(t_n)  + 
\lim_{\delta\to0^{+}} \arg L\(\tfrac{1}{2}+\delta+it_n,\chi\) 
- \dfrac{\arg G\(\chi\)}{2} = \(n - \dfrac{1}{2} - \dfrac{a}{4} \)\pi.
\eeq

\subsection{Asymptotic equation for the $n$-th zero}
\label{sec:dirichlet_asymptotic}

From Stirling's formula we have the following 
asymptotic form for $t \to \pm \infty$:
\beq
\label{approxRS}
\vartheta_{k,a} (t) = \mathrm{sgn}(t)
\left[ \dfrac{|t|}{2} \log\(  \dfrac{k |t|}{2 \pi e} \)  +  
\dfrac{2a-1}{8} \pi
+O(1/t) \right].
\eeq
The first order approximation of \eqref{exact}, i.e.
neglecting $O(1/t)$ terms, is therefore given by
\beq\label{asymptotic}
\begin{split}
\sig_n  \dfrac{|t_n|}{2\pi}\log\(\dfrac{k \, |t_n|}{2\pi e }\)
&+ \dfrac{1}{\pi}\lim_{\delta\to 0^{+}}
\arg L\(\tfrac{1}{2}+ \delta + i \sig_n |t_n|, \chi\) \\
&- \dfrac{1}{2\pi}\arg G\(\chi\)
 = n + \dfrac{\sig_n -4 - 2a(1+\sig_n)}{8},
\end{split}
\eeq
where $\sig_n = 1$ if $n>0$ and $\sig_n = -1$ if $n\le 0$. For $n>0$
we have $t_n = |t_n|$ and for $n\le 0$ we have $t_n = - |t_n|$.

\subsection{Counting formulas}
\label{sec:dirichlet_counting}

Let us define $N^+_0 (T,\chi)$ as  the number of zeros on the critical 
line with $0< \Im( \rho) < T$ and  $N^-_0 (T,\chi)$ as the number of zeros 
with $-T < \Im( \rho)  < 0$. As explained before,
$N^+_0(T,\chi) \neq N^-_0 (T,\chi)$
if the characters are complex numbers, since  the zeros are not 
symmetrically distributed between the upper and lower half of 
the critical line. 

\begin{figure}
\centering
\includegraphics[width=0.6\linewidth]{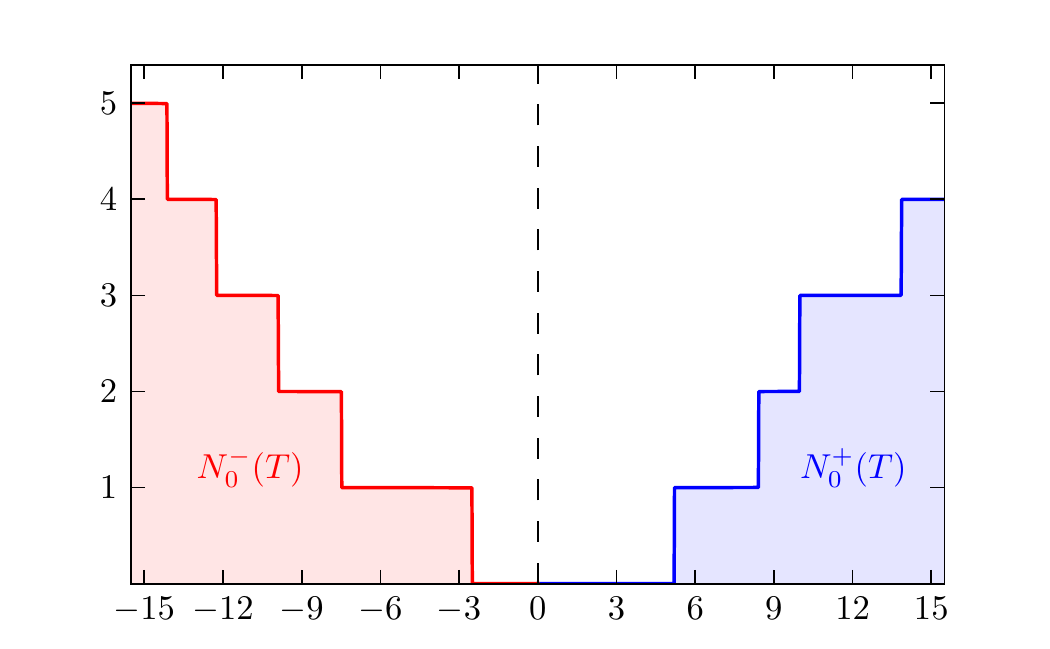}
\caption{Exact counting formulae \eqref{count_exact1}  
and \eqref{count_exact2}.
Note that they are not symmetric with respect to the origin, since the
$L$-zeros for complex $\chi$ are not complex conjugates.
We used $\chi=\chi_{7,2}$ shown in \eqref{char72}.}
\label{fig:counting}
\end{figure}

The counting formula $N^+_0(T,\chi)$  is obtained 
from equation \eqref{exact} by replacing $t_n \to T$ 
and $n \to N^+_0 + \tfrac{1}{2} $, 
therefore
\beq\label{count_exact1}
N^+_0(T,\chi) = \dfrac{1}{\pi}\vartheta_{k,a}(T) + 
\dfrac{1}{\pi}\arg L\(\tfrac{1}{2}+ i T, \chi\)
- \dfrac{1}{2\pi}\arg G\(\chi\) + \dfrac{a}{4}.
\eeq
The passage from \eqref{exact} to \eqref{count_exact1} is 
justified under the assumptions
already discussed in connection with 
\eqref{counting2} and \eqref{counting2_exact}, i.e. 
assuming that \eqref{exact} has a unique solution for every $n$. 
As explained above for the Riemann $\zeta$ case,  this
is equivalent to assume the existence of the 
$\lim_{\delta\to0^+} \arg L\(\tfrac{1}{2}+\delta+it, \chi\)$ for
every $t$.
Analogously, the counting formula on the lower half line is
given by
\beq\label{count_exact2}
N_0^-(T,\chi) = \dfrac{1}{\pi}\vartheta_{k,a}(T) - 
\dfrac{1}{\pi}\arg L \(\tfrac{1}{2} - i T, \chi\)
+\dfrac{1}{2\pi}\arg G(\chi) - \dfrac{a}{4}.
\eeq
Note that in \eqref{count_exact1} and \eqref{count_exact2} $T$ is positive.
Both cases  are plotted in Figure~\ref{fig:counting} for the character 
$\chi_{7,2}$ shown in \eqref{char72}.
One can notice that they are precisely staircase functions, jumping by one
at each zero. Note also that the functions
are not symmetric about  the origin.

From \eqref{approxRS} we also have the first order approximation
for $T \to \infty$,
\beq\label{count_approx}
N^+_0(T,\chi) = \dfrac{T}{2\pi}\log\(\dfrac{k \, T}{2\pi e }\) + 
\dfrac{1}{\pi}\arg L\(\tfrac{1}{2}+i T, \chi\) 
-\dfrac{\arg G}{2\pi} - \dfrac{1}{8} + \dfrac{a}{2} .
\eeq
Analogously, for the lower half line we  have
\beq\label{count_approx2}
N^-_0(T,\chi) = \dfrac{T}{2\pi}\log\(\dfrac{k \, T}{2\pi e }\) - 
\dfrac{1}{\pi}\arg L\(\tfrac{1}{2}-i T, \chi\)
+\dfrac{\arg G}{2\pi} - \dfrac{1}{8}.
\eeq
As  in \eqref{exact} again we are omitting $n_0$ since in the cases 
below $n_0=0$, but for other cases one may need to include $\pm n_0$ on
the right hand side of $N_0^{\pm}$, respectively.

It is known that the number of zeros on the \emph{entire critical strip} 
up to height $T$, i.e. in the region
$\{0 < \sigma < 1,\, 0 < t < T\}$, is given by \cite{Montgomery}
\beq
\label{count_strip}
N^+(T,\chi) = \dfrac{1}{\pi} \vartheta_{k,a}\(T\)
+\dfrac{1}{\pi}\arg L\(\tfrac{1}{2}+iT,\chi\) -
\dfrac{1}{\pi}\arg L\(\tfrac{1}{2},\chi\).
\eeq
From Stirling's approximation and noticing that 
$2a-1=-\chi(-1)$, then for $T\to\infty $ we obtain the asymptotic
approximation \cite{Selberg1,Montgomery}
\beq\label{selberg_counting}
\begin{split}
N^+(T,\chi) &= \dfrac{T}{2\pi}\log\(\dfrac{k \, T}{2\pi e} \)
+ \dfrac{1}{\pi}\arg L\(\tfrac{1}{2} + i T,\chi\) \\
& - \dfrac{1}{\pi}\arg L\(\tfrac{1}{2},\chi\) 
- \dfrac{\chi(-1)}{8} + O(1/T).
\end{split}
\eeq
Both formulas \eqref{count_strip} and \eqref{selberg_counting} are
exactly the same as \eqref{count_exact1} and \eqref{count_approx}, 
respectively. This can be seen as follows. 
From \eqref{FExi} we conclude that $\xi$ is real on the critical line.
Thus $\arg\xi\(\tfrac{1}{2}\)=0=
-\tfrac{1}{2}\arg G\(\chi\)+\arg L\(\tfrac{1}{2},\chi\)+\tfrac{\pi a}{4}$.
Then, replacing $\arg G\(\chi\)$ in \eqref{exact} we obtain
\beq\label{exact_v2}
\vartheta_{k,a}\(t_n\) + 
\lim_{\delta\to0^+}\arg L\(\tfrac{1}{2}+\delta+it_n, \chi\)
-\arg L\(\tfrac{1}{2},\chi\) = \(n -\tfrac{1}{2}\)\pi.
\eeq
Replacing $t_n\to T$ and $n\to N_0^+ + \tfrac{1}{2}$ in 
\eqref{exact_v2} we have precisely the expression
\eqref{count_strip}, and also \eqref{selberg_counting} for $T\to \infty$.
Then we conclude that $N_0^+(T,\chi) = N^+(T,\chi)$ exactly.
From \eqref{xiconju} we see that negative zeros for character
$\chi$ correspond to positive zeros for character $\overline{\chi}$.
Thus for $-T < \Im(\rho) < 0$ the counting on the critical strip
also coincides  with the counting on the critical line, since
$N_0^-(T,\chi) = N_0^+(T,\overline{\chi})$ and 
$N^-(T,\chi) = N^+(T,\overline{\chi})$.
Therefore, the number of zeros on the entire critical strip
is the same as the number of zeros on the critical line obtained
as solutions of \eqref{exact}, under the assumption that
\eqref{exact} has a unique solution for every $n$.  This is
equivalent to stating  that  
$\lim_{\delta\to0^+}\arg L\(\tfrac{1}{2}+\delta+it,\chi\)$ exists for
every $t$. This will be further exemplified in Section~\ref{sec:davenport}.

\section{Zeros of $L$-functions based on modular forms}
\label{sec:modular}

Let us  generalize the previous results to $L$-functions based on level 
one modular forms.   We first recall some basic definitions and properties.  
The modular group can be represented by the set
of $2 \times 2$ integer matrices
\beq \label{SL2Z}
\SL =  \left\{ 
A = \begin{pmatrix}[0.7]
a & b \\
c & d
\end{pmatrix} \, \bigg\vert \, \,  a,b,c,d \in \mathbb{Z}, \ \det A = 1
\right\},
\eeq
provided each matrix $A$ is identified with $-A$, i.e. $\pm A$ are regarded
as the same transformation. Thus for $\tau$ in the upper half 
complex plane, it transforms as 
$\tau  \mapsto A \tau = \tfrac{a \tau + b}{c\tau + d}$ 
under the action of the modular group.
A modular form $f$ of weight $k$ is a function that is analytic in the 
upper half complex plane which satisfies the functional relation 
\cite{Apostolmodular} 
\beq \label{fform}
f\(\dfrac{a\tau  + b}{c\tau + d} \) =  \(c\tau +d\)^k   f(\tau).
\eeq
If the above equation is satisfied for all of $\SL$,  
then $f$ is referred to as being of level one.
It is possible to define higher level modular forms  which satisfy 
the above equation for a subgroup of $\SL$.   Since our results are easily
generalized to the higher level case,  henceforth we will  only consider 
level one forms.     

For the $\SL$ element 
$\renewcommand\arraycolsep{2pt}
\begin{pmatrix}[0.6] 1 & 1 \\ 0 & 1 \end{pmatrix}$, 
the above equation \eqref{fform} implies  
the periodicity $f(\tau ) = f(\tau+1)$, 
thus it has a Fourier series
\beq \label{Fourier}
f(\tau ) = \sum_{n=0}^\infty  a_f (n) \, q^n  ,  
\qquad q \equiv  e^{2 \pi i \tau}.
\eeq
If $a_f(0) = 0$ then $f$ is called a cusp form. 

From the Fourier coefficients, one can define the Dirichlet series
\beq \label{Lmod}
L_f \(s\) = \sum_{n=1}^\infty \dfrac{a_f \(n\)}{n^s}.
\eeq
The functional equation for $L_f\(s\)$ relates it to 
$L_f\(k-s\)$,  so that the critical line is $\Re (s) = \tfrac{k}{2}$, 
where $k \ge 4$ is an even integer.   One can always shift the critical line
to $\tfrac{1}{2}$ by  replacing $a_f (n) \to a_f (n)/n^{(k-1)/2}$,  
however we will not do this here.   
Let us define 
\beq \label{ximod}
\Lambda_f(s) \equiv \(2\pi\)^{-s} \, \Gamma \( s \) \, L_f (s).
\eeq
Then the functional equation is given by \cite{Apostolmodular}
\beq \label{LFE}
\Lambda_f(s) = (-1)^{k/2} \Lambda_f(k-s).
\eeq

There are only two cases to consider since $\tfrac{k}{2}$ can be
an even or an odd  integer.  As in \eqref{xi} we can absorb
the extra minus sign  factor for the odd case. 
Thus we define $\xi_f(s) \equiv \Lambda_f(s)$ for
$\tfrac{k}{2}$ even, and then $\xi_f(s) = \xi_f(k-s)$. 
For $\tfrac{k}{2}$ odd 
we define $\xi_f(s) \equiv e^{-i \pi/2}\Lambda_f(s)$ 
implying $\xi_f(s) = \overline{\xi_f}(k-s)$. 
Representing $\xi_f(s) = |\xi_f| \, e^{i \theta }$ where 
$s=\sigma+i t$, we follow exactly the same steps as in the previous sections.
From the particular solution \eqref{particular} we conclude
that there are infinite zeros
on the critical line $\Re(\rho) = \tfrac{k}{2}$ determined by
$\lim_{\delta\to0^+}\theta\(\tfrac{k}{2}+\delta,t\)=
\(n-\tfrac{1}{2}\)\pi$.
Therefore, these zeros are given in the form
$\rho_n = \tfrac{k}{2} + it_n$, where $t_n$ is the solution
of the equation
\beq \label{exactmod}
\vartheta_k (t_n) + \lim_{\delta \to 0^+} 
\arg L_f \( \tfrac{k}{2} + \delta + i t_n \) =  
\( n - \dfrac{1 + (-1)^{k/2}}{4} \) \pi  
\eeq
where $n=1,2,\dotsc$ and we have defined
\beq \label{RSL}
\vartheta_k (t)  \equiv \arg \Gamma\( \tfrac{k}{2} + i t \) - t \log 2\pi.
\eeq
This implies that the number of solutions of \eqref{exactmod} with 
$0< t < T$ is given by
\beq \label{NTmod}
N_0\(T\) =  \inv{\pi} \vartheta_k (T) + 
\inv{\pi} \arg L_f \( \tfrac{ k}{2} + iT \) - \dfrac{1 - (-1)^{k/2}}{4}.
\eeq
In the limit of large $t_n$, neglecting terms of $O(1/t)$,
the equation \eqref{exactmod} becomes 
\beq \label{asymod}
t_n \log \(  \frac{t_n}{2 \pi e} \)  + 
\lim_{\delta \to 0^+}  \arg L_f \( \tfrac{k}{2} +\delta + i t_n\) 
= \( n - \dfrac{k+ (-1)^{k/2}}{4} \) \pi .
\eeq

\section{Approximate zeros in terms of the Lambert $W$-function}
\label{sec:lambert}

\subsection{Explicit formula}

We now show that it is possible to obtain an approximate solution to
the previous transcendental equations with an explicit formula.  
In this approximation,  there is indeed a unique solution to the equation
for every $n$.  
Let us introduce the Lambert $W$-function \cite{Corless}, which is 
defined for any complex number $z$ through the equation
\beq
\label{lambert_definition}
W(z) e^{W(z)} = z.
\eeq
The multi-valued $W$-function cannot be expressed in terms of other
known elementary functions. 
If we restrict attention to real-valued $W(x)$ there are two branches. 
The principal branch occurs when $W(x) \ge -1$ and is denoted by $W_0$, 
or simply $W$ for short, and its domain is $x \ge -e^{-1}$. 
The secondary branch, denoted by $W_{-1}$, satisfies $W_{-1}(x) \le -1$ 
for $-e^{-1} \le x < 0$.
Since we are interested only in positive real-valued solutions, we just 
need the principal 
branch where $W$ is single-valued.

Let us  start with the zeros of the $\zeta$-function, described
by equation \eqref{FinalTranscendence}. Consider its leading order 
approximation, or equivalently its average since
$\langle \arg\zeta\(\tfrac{1}{2}+iy\) \rangle = 0$. Then we have
the transcendental equation
\beq \label{ApproxTranscendence}
\dfrac{\tilde{t}_n}{2\pi}\log\(\dfrac{\tilde{t}_n}{2\pi e}\) = 
n - \dfrac{11}{8}.
\eeq
Through the transformation 
$\tilde{t}_n = 2\pi\(n-\tfrac{11}{8}\)x_n^{-1}$, this
equation can be written as 
$x_n e^{x_n} = e^{-1} \(n-\tfrac{11}{8}\)$. Comparing with 
\eqref{lambert_definition} we thus we obtain 
\beq \label{Lambert}
\tilde{t}_n = 
\dfrac{2\pi\(n-\tfrac{11}{8}\)}{W\[e^{-1}\(n-\tfrac{11}{8}\)\]}
\eeq
where $n=1,2,\dotsc$.

Although the inversion from \eqref{ApproxTranscendence} to \eqref{Lambert}
is rather simple, it is very  convenient since  it is indeed an 
explicit formula depending only on $n$,  
and $W$ is included in most numerical packages.
It gives an approximate solution for the ordinates of the Riemann zeros
in closed form. The values computed from \eqref{Lambert} are much closer 
to the Riemann zeros than Gram points,  
and one does not have to deal with violations of 
Gram's law; see Remark~\ref{gram}.     

Analogously, for Dirichlet $L$-functions, after neglecting
the $\arg L$ term, the equation \eqref{asymptotic} yields a transcendental
equation which can be written as $x_n e^{x_n} = k A_{n} e^{-1}$ through
the transformation $|t_n| = 2 \pi A_{n} x_n^{-1}$, where
\beq
A_{n}\(\chi\) = \sig_n \(n + \dfrac{\arg G(\chi)}{2\pi}\)
+\dfrac{1 - 4\sig_n - 2a\(\sig_n + 1\)}{8}.
\eeq
Thus the approximate solution is explicitly given by
\beq\label{approx}
\tilde{t}_n = \dfrac{2\pi \sig_n A_{n}\(\chi\)}{
W\[ k \, e^{-1} A_{n}\(\chi\) \]} 
\eeq
where $n=0,\pm1,\pm2,\dotsc$.
In the above formula $n=1,2,\dotsc$ correspond to positive $t_n$
solutions,  while $n=0,-1,\dotsc$ correspond to negative $t_n$ solutions.
Contrary to the $\zeta$-function, in general, the zeros are not
conjugate related along the critical line.

In the same way, ignoring the small $\arg L_f$ term in \eqref{asymod},
the approximate solution for the imaginary part of the zeros of $L$-functions 
based on level one modular forms is given by
\beq \label{yLambertmod}
\tilde{t}_n = \dfrac{A_n \pi }{W\[ (2e)^{-1}A_n \]} , \qquad
A_n = n - \dfrac{k+(-1)^{k/2}}{4} ,
\qquad n=1,2,\dotsc.
\eeq

\subsection{Further remarks}
\label{sec:lambert_remarks}

Let us focus on the approximation \eqref{Lambert} 
regarding zeros of the $\zeta$-function.
Obviously the same arguments apply to the zeros of the other classes of 
functions based on formulas \eqref{approx} and  \eqref{yLambertmod}.

\begin{remark}
The estimates given by \eqref{Lambert} can be calculated 
for arbitrarily large $n$, since $W$ is a standard elementary function.  
Of course the  $\tilde{t}_n$ are not as accurate as the solutions 
$t_n$ including the $\arg \zeta$ term, as we will see in 
Section~\ref{sec:numerical}.
Nevertheless, it is indeed a good estimate, especially if one considers 
very high zeros  where  traditional methods have 
not previously estimated  such high values.
For instance, formula \eqref{Lambert} can easily estimate the 
zeros shown in Table~\ref{highn} (Appendix~\ref{sec:tables_zeta}), 
and much higher if desirable.
The numbers in this table are accurate approximations to the $n$-th zero 
to the number of digits shown,  which is approximately the number of 
digits in the integer part.
For instance,  the approximation to the $10^{100}$ zero is
correct to $100$ digits.   
With Mathematica we easily calculated the first million digits
of the $10^{10^6}$ zero.   
\end{remark}

\begin{remark}Using the asymptotic behaviour 
$W(x) \approx \log x$ for large $x$,
the $n$-th zero is approximately given by 
$\tilde{t}_n \approx  2 \pi n/ \log n$,
as already known \cite{Titchmarsh}.
The distance between consecutive ordinates  
is then approximately equal to 
$\tilde{t}_{n+1}-\tilde{t}_n \approx 2\pi/\log n$, which
tends to zero when $n\to \infty$.
\end{remark}

\begin{remark}\label{gram}
The solutions \eqref{Lambert} are reminiscent of the so-called Gram 
points $g_n$,  which are solutions to 
$\vartheta (g_n) = n\pi$ where $\vartheta$ is given by \eqref{riemann_siegel}. 
Gram's law is the tendency for Riemann zeros to lie between 
consecutive Gram points, but it is known to fail for about $\tfrac{1}{4}$ 
of all 
Gram intervals. Our $\tilde{t}_n$ are intrinsically different from 
Gram points. It is an approximate solution for the ordinate of
the zero itself. In particular,  the  Gram point $g_0 = 17.8455$ 
is the closest to 
the first Riemann zero, whereas $\tilde{t}_1 = 14.52$ is 
already much closer to the true zero which is $t_1 \approx 14.1347$.
The traditional method to compute the zeros is based on 
the Riemann-Siegel formula 
$\zeta\(\tfrac{1}{2}+it\) = Z(t)\[\cos\vartheta(t) - i \sin\vartheta(t)\]$,
and the empirical observation that the real part of this equation 
is almost always positive, except when Gram's law fails, and $Z(t)$ has 
the opposite
sign of $\sin\vartheta$. Since $Z(t)$ and $\zeta\(\tfrac{1}{2}+it\)$ have
the same zeros, one looks for the zeros of $Z(t)$ between two
Gram points, as long as Gram's law holds $(-1)^nZ\(g_n\)>0$. 
To verify the RH numerically, the counting formula \eqref{backlund} must 
also be used to assure that the number of zeros on the critical line 
coincide with the number of zeros on the strip. 
The detailed procedure is throughly 
explained in \cite{Edwards,Titchmarsh}.
Based on this method, amazingly accurate solutions and high zeros on 
the critical line were computed 
\cite{Gourdon,Odlyzko,OdlyzkoSchonhage,Odlyzko2}.
Nevertheless, our proposal is \emph{fundamentally} different.
We claim that \eqref{exact_eq2} is the equation that determines
the Riemann zeros on the critical line.
Then, one just needs to find its solution for a given
$n$. We will compute the Riemann zeros  in this way in the next 
section, just by solving the equation \eqref{exact_eq2} numerically, starting
from the approximation given by the explicit formula \eqref{Lambert},
without using Gram points nor the Riemann-Siegel $Z$ function.
Let us emphasize that our goal is not to provide a more efficient algorithm
to compute the zeros \cite{OdlyzkoSchonhage},  although the method 
described here may very well be, but to justify the 
validity of equation \eqref{exact_eq2}.
\end{remark}

\section{A counterexample: the Davenport-Heilbronn function}
\label{sec:davenport}

In this section 
we consider a function that is known to violate the RH, and this serves to 
sharpen our understanding of our previous analysis.  
In this example,  one  can clearly see how the corresponding 
transcendental equation does  not 
have a unique solution for every $n$.

The \emph{Davenport-Heilbronn} function is  defined by
\beq\label{davenport_def}
\CD(s) \equiv \dfrac{(1-i \kappa)}{2} L\(s,\chi_{5,2}\) +
\dfrac{(1+i \kappa)}{2} L\(s,\overline{\chi}_{5,2}\)
\eeq
with 
\beq
\kappa =\dfrac{\sqrt{ 10-2\sqrt{5}}-2}{\sqrt{5} - 1}.
\eeq
Above the Dirichlet character is the following:
\beq\label{char52}
\begin{tabular}{@{}c|ccccccc@{}}
$n$             & $1$ & $2$ & $3$ & $4$ & $5$ \\
\midrule[0.3pt]
$\chi_{5,2}(n)$ &
$1$ & $i$ & $-i$ & $-1$ & $0$
\end{tabular}
\eeq
where $\chi_{5,2}(-1) = -1$ thus $a=1$. The function \eqref{davenport_def}
satisfies the functional equation
\beq\label{davenport_func_eq}
\xi(s) = \xi(1-s), \qquad
\xi(s) \equiv \(\dfrac{\pi}{5}\)^{-s/2}\, \Gamma\(\dfrac{1+s}{2}\) \CD(s).
\eeq

The function \eqref{davenport_def} has almost all the same
properties of $\zeta$, such as a functional equation, except that it has
no Euler Product Formula. It is well known that it has zeros in 
the region $\Re\(s\) > 1$,   which is essentially a  consequence that it has
no Euler product.   It also has zeros
in the critical strip  $0 \le \Re(s) \le 1$, where infinitely many of them
lie on the critical line $\Re(s) = \tfrac{1}{2}$, however,  it also has
zeros off of the critical line, thus violating the RH.
For a detailed study of this function and numerical computation
of its zeros see \cite{BombieriGosh}.

Repeating the analysis of the previous sections for zeros on 
the critical line, we obtain the following
transcendental equation:
\beq\label{davenport_transcendental}
\dfrac{1}{\pi}\vartheta_{5,1}(t_n) + 
\dfrac{1}{\pi}\lim_{\delta\to0^+} \arg \CD\(\tfrac{1}{2}+\delta+it_n\) 
+\dfrac{1}{2} = n
\eeq
where $\vartheta_{5,1}$ is defined in \eqref{RSgen}. The approximate
solution is explicitly given by
\beq\label{davenport_lambert}
\tilde{t}_n = \dfrac{2\pi\(n-\tfrac{5}{8}\)}{W\[ 5e^{-1}(n-\tfrac{5}{8}) \]}
\eeq
for $n=1,2,\dotsc$.
From \eqref{davenport_transcendental} and
\eqref{davenport_lambert} it is possible to compute
zeros on the critical line. 
Moreover, zeros off of the critical line
satisfy the general solution \eqref{gen_sol} \cite{Lectures}. This
shows that $B=0$, with $B$ defined in \eqref{Bdef}, captures
all the zeros.

Since \eqref{davenport_transcendental} only captures zeros at 
$\sigma=\tfrac{1}{2}$, what happens if there are zeros off of 
the critical line?
Consider a simple zero  
denoted by $\rho_{\bullet} = \sigma_{\bullet} + i \tzero$
where $0< \sigma_{\bullet} < 1$ and $\sigma_{\bullet} \ne \tfrac{1}{2}$.
Due to the functional equation \eqref{davenport_func_eq} there is also a zero 
at $1-\overline{\rho}_{\bullet} = 1-\sigma_{\bullet} + i \tzero$.
Let 
\beq\label{Sdav}
S_{\CD}(t) = \dfrac{1}{\pi}\lim_{\delta\to0^+}
\arg\CD\(\tfrac{1}{2}+\delta+it\). 
\eeq
From its role 
in the counting formula over the entire critical strip,  
one knows that when $t$ varies
across $\tzero$ then $S_{\CD}(t)$ must jump by two, i.e.
we must have
$\Delta S_{\CD}(\tzero) \equiv
S_{\CD}(\tzero+\epsilon) - S_{\CD}(\tzero-\epsilon) = 2$.
This implies that $S_\CD(t)$ changes branch 
around $\tzero$ in such a way that
it cannot be smoothed out by the $\delta\to0^+$ limit. In other words,
the limit \eqref{Sdav} does
not exist close to $\tzero$. 
Therefore, \eqref{davenport_transcendental} will not have 
a solution around $\tzero$ for a given $n$.
If instead of a simple zero we have a zero with multiplicity
$m \ge 2$, then $\Delta S_{\CD}(\tzero) = 2 m$, changing branch
even more drastically.
The same situation also happens  if there are zeros 
with multiplicity $m\ge 2$ on the critical line, where we
would have $\Delta S_{\CD} = m$.

\begin{figure}
\centering
\begin{minipage}{.5\textwidth}
  \centering
  \includegraphics[width=1\linewidth]{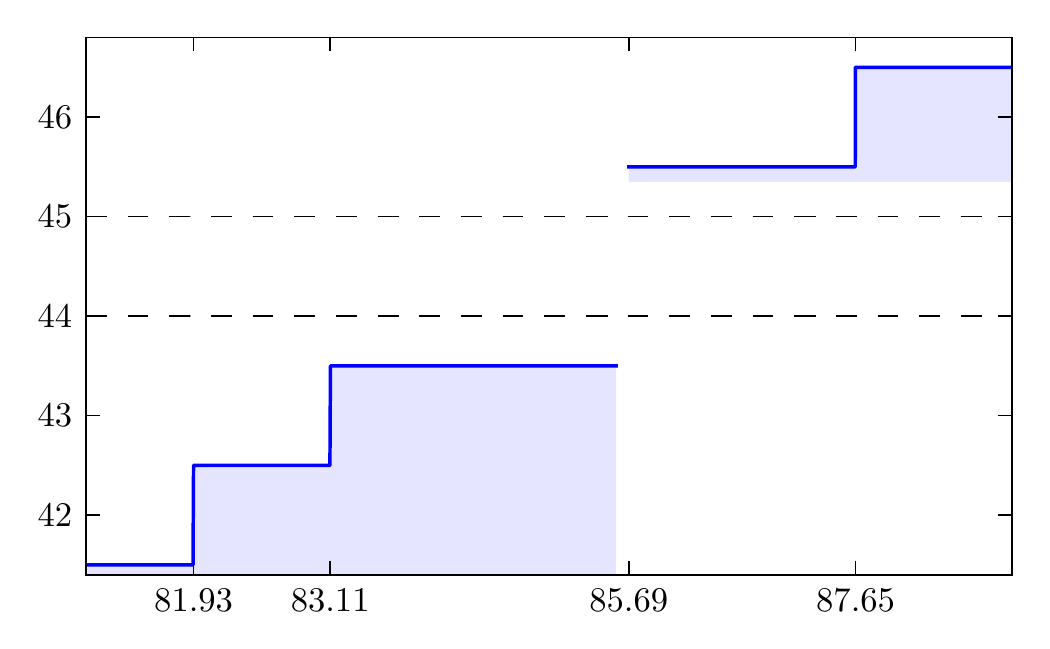}\\[-1em]
  \small{(a)}
\end{minipage}%
\begin{minipage}{.5\textwidth}
  \centering
  \includegraphics[width=1\linewidth]{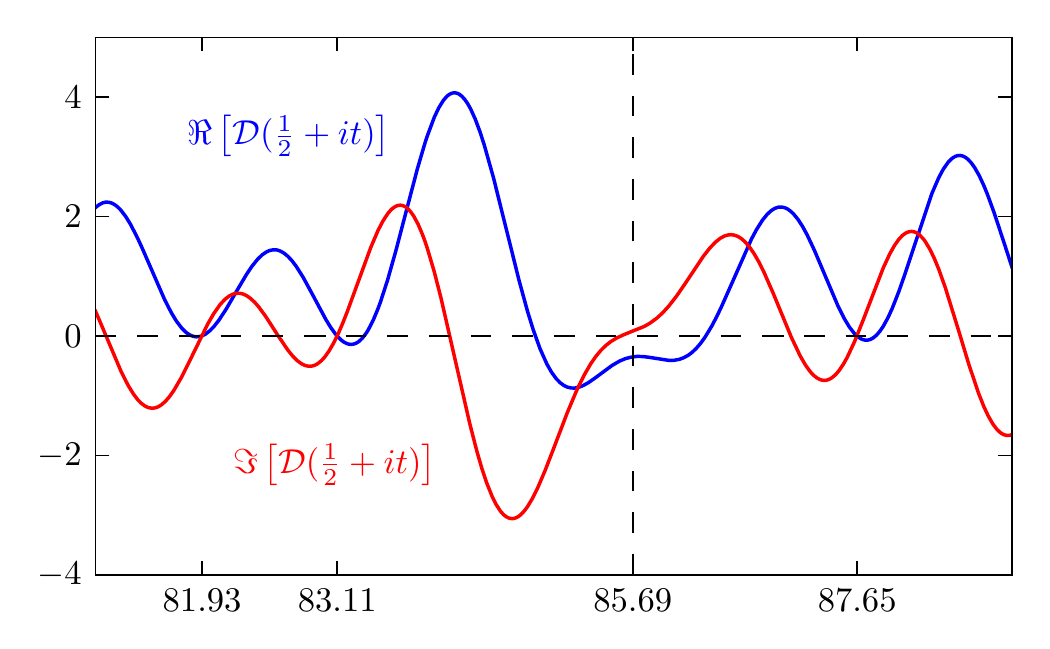}\\[-1em]
  \small{(b)}
\end{minipage}
\caption{
(a) Left hand side of \eqref{davenport_transcendental}
against $t$. Note the discontinuity
at the point $\tzero\approx 85.6993$ corresponding to $n=44$ and $n=45$, where
\eqref{davenport_transcendental} has no solution.
(b) We plot $\Re \[ \CD(1/2+it) \]$ (blue line) 
and $\Im \[ \CD(1/2+it) \]$ (red line) against $t$. 
Observe that $\Re\[ \CD(1/2+i t ) \] < 0$ when 
$\Im\[\CD(1/2+i t)\] \to 0$, for $t \to \tzero$, 
signaling the change of branch of $\arg\CD(1/2+i t)$.
}
\label{fig:davenport}
\end{figure}

In the case of the  function \eqref{davenport_def}, 
the first zero off of the critical
line occurs at $\sigma_{\bullet} \approx 0.8085$ and $\tzero \approx 85.6993$.
In Figure~\ref{fig:davenport}a we plot the left hand 
side of \eqref{davenport_transcendental} against $t$, and 
one can clearly see the above mentioned situation, 
namely that \eqref{davenport_transcendental} is not defined
at $\tzero$ and there is no solution for $n=44$ and $n=45$. The change
of branch close to $\tzero$ can be seen from Figure~\ref{fig:davenport}b.
Therefore, denoting $N_0(T)$ the number of solutions of
\eqref{davenport_transcendental} up to height $T$, we clearly
have $N_0(T) < N(T)$, where $N(T)$ is the number of zeros in
the entire critical strip. For a more detailed illustration of these
facts we refer the reader to \cite{Lectures}.

For simple zeros on the critical line the limit \eqref{Sdav} exists,
but for zeros off of the critical line,   it does not, since $S_\CD(t)$ has
to jump at least by two and the change of branch does not allow us
to smooth the function.

\section{Numerical analysis: $\zeta$-function}
\label{sec:numerical}

\subsection{The importance of $\arg \zeta$} 

Instead of solving the exact equation \eqref{exact_eq2} we will 
initially consider its first order approximation, 
which is equation \eqref{FinalTranscendence}. 
As we will see, this approximation already yields surprisingly accurate 
values for the Riemann zeros.

\begin{figure}
\centering
\begin{minipage}{.5\textwidth}
  \centering
  \includegraphics[width=1\linewidth]{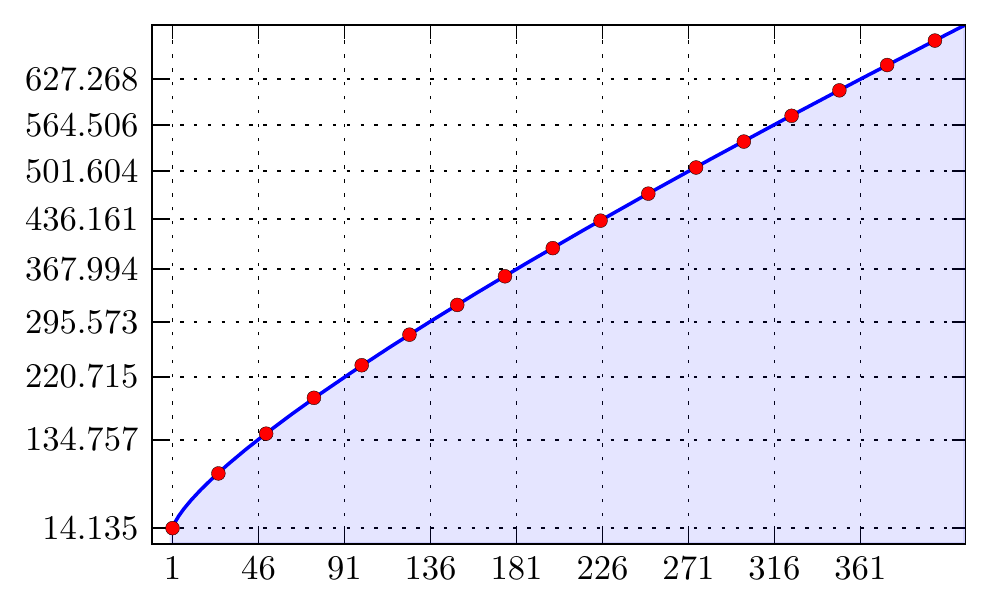}\\[-1em]
  \small{(a)}
\end{minipage}%
\begin{minipage}{.5\textwidth}
  \centering
  \includegraphics[width=1\linewidth]{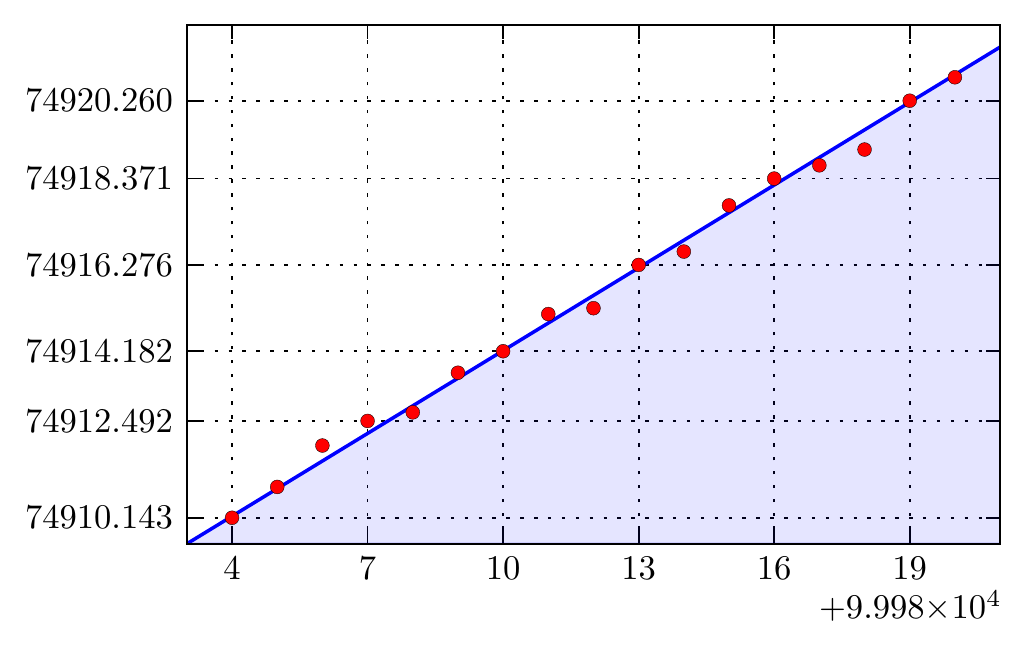}\\[-1em]
  \small{(b)}
\end{minipage}
\caption{Comparison of the prediction of \eqref{Lambert} (blue line) and 
\eqref{FinalTranscendence} (red dots). We are plotting $t_n$ against $n$.
(a) $n \in [1,\dotsc,400]$. Note how the solutions are close at first
sight. (b) If we focus on a small range we can see how
the solutions of \eqref{FinalTranscendence} oscillate around the line  
\eqref{Lambert} due to the fluctuating term $\arg \zeta$. Here
$n \in [99984,\dotsc,10^5]$.}
\label{fig:trans_zeros}
\end{figure}

Let us first consider how the approximate solution given by
\eqref{Lambert} is modified by the presence of the $\arg\zeta$ term
in \eqref{FinalTranscendence}.  Numerically, we compute $\arg\zeta$ taking
its principal value.   The fact that we get very accurate zeros up to 
the billionth zero implies that up to this $t$,  $\arg \zeta$ near a zero is 
always on the principal branch.   
As already discussed in Remark~\ref{arg_term}, the function 
$\arg\zeta\(\tfrac{1}{2} + i t\)$ oscillates around its average, which
is zero,   
as shown in Fi\-gu\-re~\ref{fig:arg_counting}a.
At a Riemann zero it can be defined by the limit \eqref{deltadef} 
which is generally not zero. 
The $\arg \zeta$ term plays an important role and 
indeed improves the estimate of the $n$-th zero.
This can be seen in Figure~\ref{fig:trans_zeros} where we compare the
estimate given by \eqref{Lambert} with the numerical solutions of
\eqref{FinalTranscendence}.

We can apply a root finder 
method in an appropriate interval, 
centered around the approximate solution $\tilde{t}_n$ given by 
formula \eqref{Lambert}.  
Some of the solutions obtained in this way are presented in 
Table~\ref{some_zeros} (Appendix~\ref{sec:tables_zeta}) 
and are accurate up to the number of decimal places
shown. We used  only Mathematica or some very simple algorithms to perform 
these  numerical computations, taken from standard open source numerical 
libraries. 

Although the equation \eqref{FinalTranscendence} was 
derived for large $n$, it is 
surprisingly accurate even for the lower zeros, as shown in
Table~\ref{lower_zeros} (Appendix~\ref{sec:tables_zeta}).   
It is actually easier to solve for low 
zeros since  $\arg \zeta $ is better behaved.    
These numbers are correct 
up to the number of digits shown, and the precision was improved simply 
by decreasing the error tolerance.

\subsection{GUE statistics}
\label{sec:gue}

The link between the Riemann zeros and random matrix theory
started with the pair correlation of zeros, proposed by 
Montgomery \cite{Montgomery},  and
the observation of Dyson \cite{Dyson} that it is the same as 
the 2-point correlation
function predicted by the Gaussian Unitary Ensemble (GUE) for large 
random  matrices.

The main purpose of this section is to test whether our 
approximation \eqref{FinalTranscendence} to the zeros 
is accurate enough to reveal this statistics.
Whereas formula \eqref{Lambert} 
is a valid estimate, 
it is not sufficiently accurate to reproduce 
the GUE statistics,  since it does not have the 
oscillatory  $\arg \zeta$ term. 
On the other hand, the solutions to equation \eqref{FinalTranscendence} 
are accurate enough,
which again  indicates the importance of $\arg \zeta$.   

Montgomery's pair correlation conjecture can be stated as follows:
\beq
\label{montgomery}
\dfrac{1}{N(T)}
\sideset{}{'}\sum_{\substack{
0\le t,t'\le T \\[0.4em]
\alpha < d(t,t') \le \beta
}} \hspace{-.7em} 
1
\, \sim \,  \int_{\alpha}^{\beta}du
\(1 - \dfrac{\sin^2\(\pi u\)}{\pi^2 u^2}\),
\eeq
where 
$d(t,t')=\tfrac{1}{2\pi}\log\(\tfrac{T}{2\pi}\)\(t-t'\)$,
$ 0 < \alpha<\beta$, 
$N(T)\sim \tfrac{T}{2\pi}\log\(\tfrac{T}{2\pi}\)$ 
according to \eqref{riemann_counting}, and the statement is valid in the limit
$T\to \infty$. The
right hand side of \eqref{montgomery} is the 2-point GUE correlation 
function. The average spacing between consecutive zeros is given by 
$\tfrac{T}{N} \sim 2\pi/\log\(\tfrac{T}{2\pi}\)\to 0$  
as $T\to \infty$. This can also be seen from \eqref{Lambert} for very
large $n$, i.e. $t_{n+1}-t_n \to 0$ as $n\to\infty$. 
Thus the distance $d(t,t')$ between zeros on the
left hand side of \eqref{montgomery} is a normalized
distance.

While \eqref{montgomery} can be applied if we start from the first
zero on the critical line, it is unable to provide a test if we are centered
around a given high zero on the line. To deal with such a situation,
Odlyzko \cite{Odlyzko2} proposed a stronger version of Montgomery's 
conjecture by taking into account the large density of zeros 
higher on the line. This is done by replacing the normalized distance in 
\eqref{montgomery} by a sum of normalized distances over
consecutive zeros in the form
\beq
d_n \equiv \dfrac{1}{2\pi}\log\(\dfrac{t_n}{2\pi}\)\(t_{n+1}-t_n\).
\eeq
Thus \eqref{montgomery} is replaced by
\beq
\label{odlyzko_pair}
\dfrac{1}{\(N-M\)\(\beta-\alpha\)}\hspace{-1em}\sideset{}{'}\sum_{\substack{
M \le m,n \le N \\[0.4em]
\alpha < \sum_{k=1}^{n}d_{m+k} \le \beta
}} \hspace{-1.4em} 1
\, = \, \dfrac{1}{\beta-\alpha}\int_{\alpha}^{\beta}du
\(1 - \dfrac{\sin^2\(\pi u\)}{\pi^2 u^2}\),
\eeq
where $M$ is the label of a given zero on the line and
$N > M$. In this sum it is also assumed that $n > m$, and we
included the correct normalization on both sides. The conjecture
\eqref{odlyzko_pair} is already  well supported by extensive 
numerical analysis \cite{Odlyzko2,Gourdon}. 

\begin{figure}
\centering
\begin{minipage}{.49\textwidth}
    \centering
    \includegraphics[width=1\linewidth]{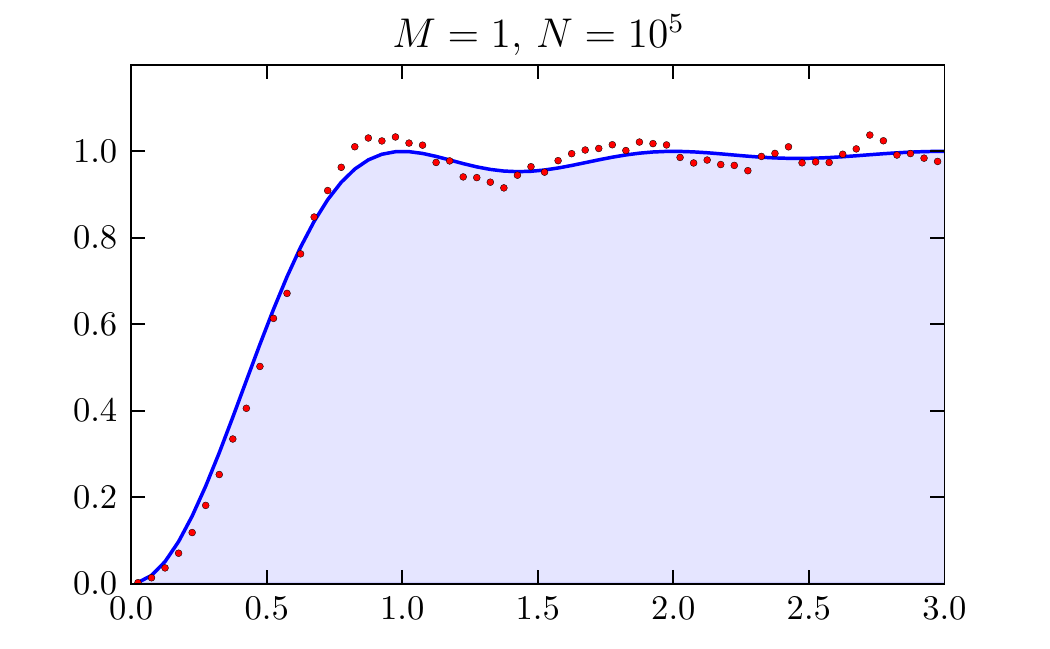}\\[-1em]
    \small{(a)}
\end{minipage}
\begin{minipage}{.49\textwidth}
    \centering
    \includegraphics[width=1\linewidth]{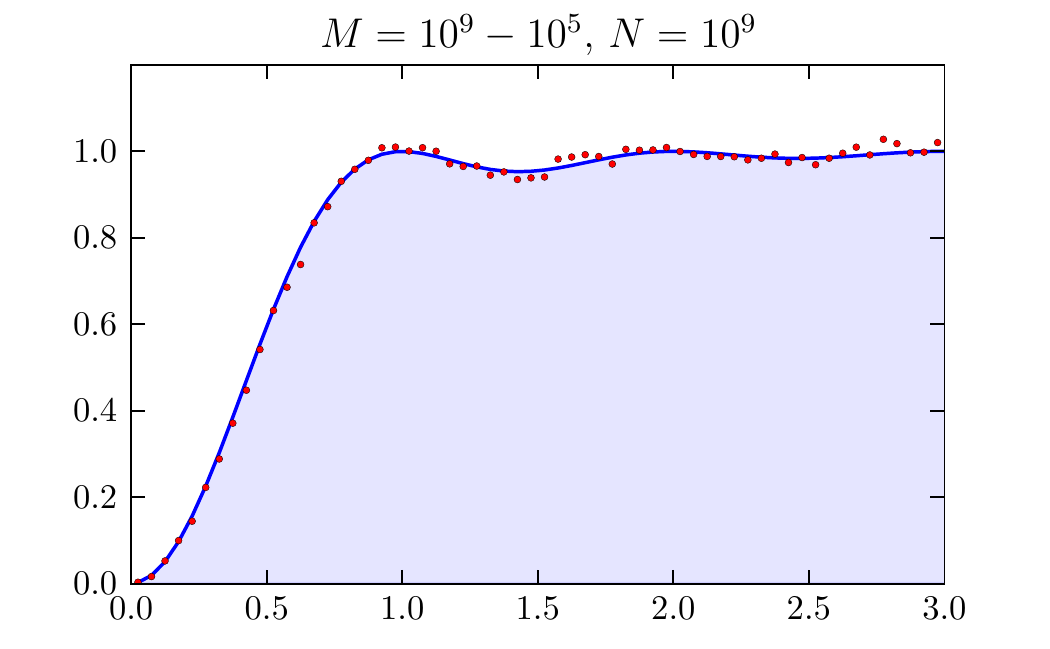}\\[-1em]
    \small{(b)}
\end{minipage}
\caption{The solid line represents the right hand side of \eqref{odlyzko_pair}
and the dots represent its left hand side, computed from
equation \eqref{FinalTranscendence}.
The parameters are $\beta = \alpha + 0.05$,
$\alpha=(0,\, 0.05,\,\dotsc,\,3)$ and the $x$-axis is given by
$x=\tfrac{1}{2}\(\alpha+\beta\)$. (a) We use the first $10^5$ zeros.
(b) The same parameters but using zeros in the middle of the
critical line; $M=10^9-10^5$ and $N=10^9$.}
\label{fig:gue}
\end{figure}

Odlyzko's conjecture \eqref{odlyzko_pair} is a very strong constraint on
the statistics of the zeros. Thus we submit the numerical solutions of 
equation \eqref{FinalTranscendence} 
to this test. In 
Figure~\ref{fig:gue}a we can see the result for $M=1$ and $N=10^{5}$, with
$\alpha$ ranging from $0\dotsc 3$ in steps of $\epsilon=0.05$, and
$\beta=\alpha+\epsilon$ for each value of $\alpha$, i.e. 
$\alpha = (0.00, \, 0.05, \, 0.10,  \dotsc,\,3.00)$ and 
$\beta = (0.05,\, 0.10, \, \dotsc, \, 3.05)$. 
We compute the left hand side of \eqref{odlyzko_pair} for each 
pair $(\alpha, \beta)$ and plot
the result against $x = \tfrac{1}{2}\(\alpha + \beta\)$.
In Figure~\ref{fig:gue}b we do the same thing but with
$M=10^9-10^5$ and $N=10^9$.
Clearly, the numerical solutions of \eqref{FinalTranscendence} reproduce
the GUE statistics. In fact, Figure~\ref{fig:gue}a is identical
to the one in \cite{Odlyzko2}. The last zeros in these ranges are
shown in Table~\ref{high_values} (Appendix~\ref{sec:tables_zeta}).

\subsection{Prime number counting function}  
\label{sec:prime}

In this section we explore whether our approximations to the Riemann zeros
are accurate enough to reconstruct the prime number counting function.     
As usual, let $\pi (x)$ denote the number of primes less than $x$.
Riemann obtained an explicit expression for $\pi (x)$ in terms of 
the non-trivial zeros of $\zeta(s)$.  
There are simpler but equivalent versions 
of the main result,  based  on the function $\psi (x) $ below.   
However, let us present the main formula for $\pi (x)$ itself
since it is historically more important.

The function $\pi (x)$ is related to another number-theoretic 
function $J(x)$, defined as 
\beq
\label{Jx}
J(x)  =  \sum_{2\leq n \leq x}    \frac{\Lambda (n)}{\log n} 
\eeq
where  $\Lambda (n)$,  the  von Mangoldt function,
is defined as $\Lambda(n) = \log p$ if $n=p^m$ for some prime $p$  
and an integer $m\ge 1$, and $\Lambda(n) = 0$ otherwise.  The two functions 
$\pi (x)$ and $J(x)$ are related by M\" obius inversion:
\beq
\label{mob1}
\pi (x) = \sum_{n\geq 1} \frac{\mu (n)}{n}  J(x^{1/n}).
\eeq
Here  $\mu (n)$ is the M\" obius function defined
as follows. $\mu(n) = 0$ if $n$
has one or more repeated prime factors, $\mu(n)=1$ if $n=1$
and $\mu(n)=(-1)^k$ if $n$ is a product of $k$ distinct primes.
The above expression is actually a finite sum,  
since for large enough $n$,  $x^{1/n} <2$ and $J=0$.

The main result of Riemann is a formula for $J(x)$, expressed 
as an infinite sum over zeros $\rho$  of the $\zeta(s)$ function:
\beq
\label{Jzeros}
J(x) =  \Li (x) - \sum_\rho  \Li\(x^\rho\)  +  
\int_x^\infty  \dfrac{dt}{\log t} ~  \inv{t \(t^2 -1\)} - \log 2,
\eeq
where $\Li (x) = \int_0^x dt /\log t $ is the log-integral 
function\footnote{Some care must be taken in  numerically 
evaluating $\Li (x^\rho )$  since $\Li$ has a branch point. 
It is more properly defined as ${\rm Ei} (\rho \log x )$ 
where ${\rm Ei} (z) = - \int_{-z}^\infty  dt \, e^{-t} /t$ is the 
exponential integral function.}.
The above sum is real because the $\rho$'s come 
in conjugate pairs.
If there are no zeros on the line $\Re(z) = 1$, then the dominant 
term is the first one in the above equation, $J(x) \approx  \Li (x)$,  
and this was used to  prove the prime number theorem by 
Hadamard and de la Vall\' ee Poussin.

The function $\psi (x)$ has the simpler form
\beq
\label{psizeros}
\psi (x)  =   \sum_{n \leq x}   \Lambda (n) = 
x - \sum_\rho\dfrac{x^\rho}{\rho} - \log (2\pi)  - \inv{2} 
\log \( 1 - \inv{x^2} \).
\eeq
In this formulation  the prime number theorem is equivalent
to $\psi (x) \approx x$.

\begin{figure}
\centering
\begin{minipage}{.5\textwidth}
    \centering
    \includegraphics[width=1\linewidth]{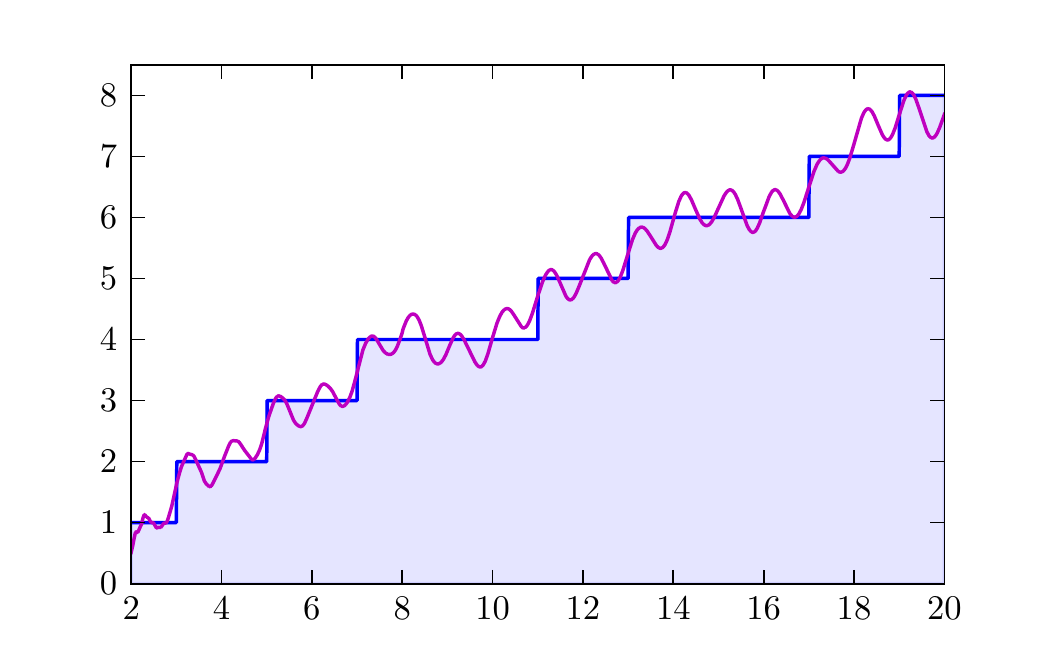}\\[-1em]
    \small{(a)}
\end{minipage}%
\begin{minipage}{.5\textwidth}
    \centering
    \includegraphics[width=1\linewidth]{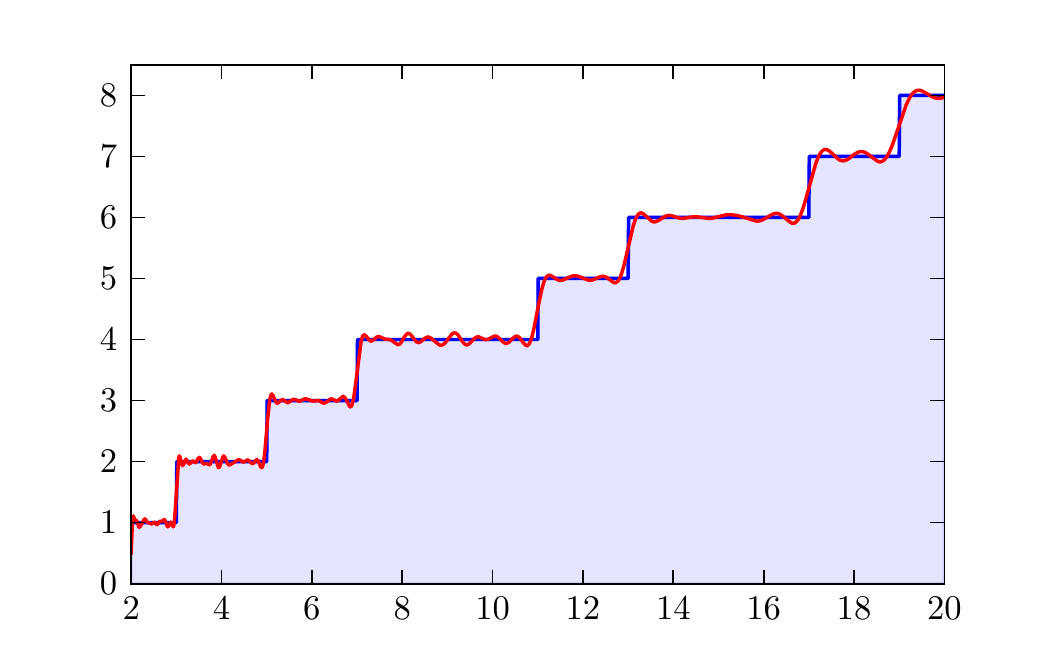}\\[-1em]
    \small{(b)}
\end{minipage}
\caption{The prime number counting function $\pi (x)$ with the 
first $50$ Riemann
zeros. (a) Zeros approximated by the formula \eqref{Lambert}.
(b) Zeros obtained from  numerical solutions to the 
equation \eqref{FinalTranscendence}.}
\label{fig:prime}
\end{figure}

In Figure~\ref{fig:prime}a we plot $\pi(x)$  from equations \eqref{mob1}
and \eqref{Jzeros}, computed with the first $50$ zeros in the 
approximation $\rho_n = \tfrac{1}{2} + i \tilde{t}_n$ given by 
\eqref{Lambert}. Figure~\ref{fig:prime}b  shows the same plot with zeros 
obtained from the numerical solutions of equation \eqref{FinalTranscendence}.
Although with the approximation $\tilde{t}_n$ the curve is trying 
to follow the steps in $\pi (x)$,  once again, one clearly sees the importance 
of the $\arg  \zeta$ term.

\subsection{Solutions to the exact equation}
\label{sec:numerical_exact}

In the previous sections we have computed numerical solutions of
\eqref{FinalTranscendence} showing that, actually, this first
order approximation to \eqref{exact_eq2} is very good and already 
captures  some  interesting properties of the Riemann zeros, 
such as the GUE statistics and the ability
to reproduce $\pi(x)$.
Nevertheless, by simply 
solving \eqref{exact_eq2} it is possible to obtain values for the zeros
as accurately  as desirable. 
The numerical procedure is performed  as follows:
\begin{enumerate}
\item \label{step1} We apply a root finder method on \eqref{exact_eq2} 
looking for the solution in a region centered around the 
number $\tilde{t}_n$ provided by \eqref{Lambert}, with a not so 
small $\delta$, for instance $\delta \sim 10^{-5}$.
\item \label{step2} We solve \eqref{exact_eq2} again but now centered around 
the solution obtained in step \ref{step1} above, and we decrease $\delta$, for 
instance $\delta \sim 10^{-8}$.
\item We repeat the procedure in step \ref{step2} above, 
decreasing $\delta$ again.
\item Through successive iterations, and decreasing $\delta$ 
each time, it is possible to obtain solutions as accurate as desirable.
In carrying this out,  
it is important to not allow $\delta$ to be exactly zero.   
\end{enumerate}
An actual implementation of the above procedure in Mathematica is shown
in Appendix~\ref{sec:mathematica},  
which we have included mainly to show its simplicity.   
The first few zeros computed in this way are shown 
in Table~\ref{lower_precise} (Appendix~\ref{sec:tables_zeta}).
Through successive iterations it is possible achieve even much higher  
accuracy than  shown in 
Table~\ref{lower_precise}.

It is known that the first zero where  Gram's law fails is for $n=126$.
Applying the same method, like for any other $n$, the solution of 
\eqref{exact_eq2} starting with the approximation \eqref{Lambert} does
not present any difficulty. We easily found the following number:
\begin{flalign*}
t_{126} = 279.229250927745189228409880451955359283492637405561293594727
\end{flalign*}
Just to illustrate, and to convince the reader,  how the solutions 
of \eqref{exact_eq2} can be made arbitrarily precise, we compute the 
zero $n=1000$ accurate up to $500$ decimal places, also using the same 
simple approach\footnote{Computing this number to $500$ digit accuracy  
took a few minutes on a standard personal laptop 
computer using Mathematica. It only takes a few seconds to 
obtain $100$ digit accuracy.}:
\begin{flalign*}
t_{1000} = 
1419.&42248094599568646598903807991681923210060106416601630\\[-.4em]
     &46908146846086764175930104179113432911792099874809842\\[-.4em]
     &32260560118741397447952650637067250834288983151845447\\[-.4em]
     &68825259311594423942519548468770816394625633238145779\\[-.4em]
     &15284185593431511879329057764279980127360524094461173\\[-.4em]
     &37041818962494747459675690479839876840142804973590017\\[-.4em]
     &35474131911629348658946395454231320810569901980719391\\[-.4em]
     &75430299848814901931936718231264204272763589114878483\\[-.4em]
     &29996467356160858436515425171824179566414953524432921\\[-.4em]
     &93649483857772253460088
\end{flalign*}
Furthermore,  one can substitute known  
precise Riemann zeros into \eqref{exact_eq2} and 
can check that the equation is identically satisfied.
These results corroborate that \eqref{exact_eq2} is an
exact equation for the Riemann zeros.

\section{Numerical analysis: $L$-functions}
\label{sec:numerical_lfunc}

We perform exactly the same numerical procedure as described
in the previous Section~\ref{sec:numerical_exact},
but now with equation \eqref{exact} and \eqref{approx} for 
Dirichlet $L$-functions, or with
\eqref{exactmod} and \eqref{yLambertmod}  for $L$-functions based 
on level one modular forms.

\subsection{Dirichlet $L$-functions}
\label{sec:numerical_dirichlet}

We will illustrate our formulas with the primitive 
characters $\chi_{7,2}$ and $\chi_{7,3}$   
since they possess the full  generality of $a=0$ and $a=1$ and 
complex components.  
There are actually $\varphi(7)=6$ distinct characters mod $7$.   

\bigskip

\noindent{\bf Example \boldmath{$\chi_{7,2}$}.}
Consider $k=7$ and $j=2$, i.e. we are computing
the Dirichlet character $\chi_{7,2}(n)$. For this case $a=1$.
Then we have the following components:
\beq\label{char72}
\def\arraystretch{1.2}
\begin{tabular}{@{}c|ccccccc@{}}
$n$             & $1$ & $2$ & $3$ & $4$ & $5$ & $6$ & $7$ \\
\hline
$\chi_{7,2}(n)$ &  
$1$ & $e^{ 2\pi i /3}$ & $e^{\pi i / 3}$ & $e^{-2\pi i / 3}$ &
$e^{-\pi i / 3}$ & $-1$ & $0$ 
\end{tabular} 
\eeq
The first few zeros, positive and negative, obtained by solving
\eqref{exact} are shown in Table~\ref{zeros_1} in
Appendix~\ref{sec:tables_dirichlet}.
The solutions shown are
easily obtained with $50$ decimal places of accuracy.

\bigskip

\noindent{\bf Example \boldmath{$\chi_{7,3}$}.} 
Consider $k=7$ and $j=3$, such that $a=0$.
In this case the components of $\chi_{7,3}(n)$ are the following:
\beq\label{char73}
\def\arraystretch{1.2}
\begin{tabular}{@{}c|ccccccc@{}}
$n$             & $1$ & $2$ & $3$ & $4$ & $5$ & $6$ & $7$ \\
\hline
$\chi_{7,3}(n)$ &  
$1$ & $e^{ -2\pi i /3}$ & $e^{2 \pi i / 3}$ & $e^{2\pi i / 3}$ &
$e^{-2 \pi i / 3}$ & $1$ & $0$ 
\end{tabular}
\eeq
The first few solutions of \eqref{exact} are shown in 
Table~\ref{zeros_2} in Appendix~\ref{sec:tables_dirichlet}
and are accurate up to $50$ decimal places.
As previously stated, the solutions to equation \eqref{exact} can be 
calculated to any desired level of accuracy. For instance,
we can easily compute the following 
number for $n=1000$, accurate to $100$ decimal places:
\begin{flalign*}
t_{1000} = 1037.&56371706920654296560046127698168717112749601359549\\[-.4em]
                &01734503731679747841764715443496546207885576444206
\end{flalign*}
We also have been able to solve the equation for high zeros to 
high accuracy, up to the millionth zero, some of which are 
listed in Table~\ref{highzeros} in Appendix~\ref{sec:tables_dirichlet}, 
and were previously unknown.

\subsection{Modular $L$-function based on Ramanujan $\tau$}
\label{sec:ramanujan}

\begin{figure}
\centering
\includegraphics[width=0.6\linewidth]{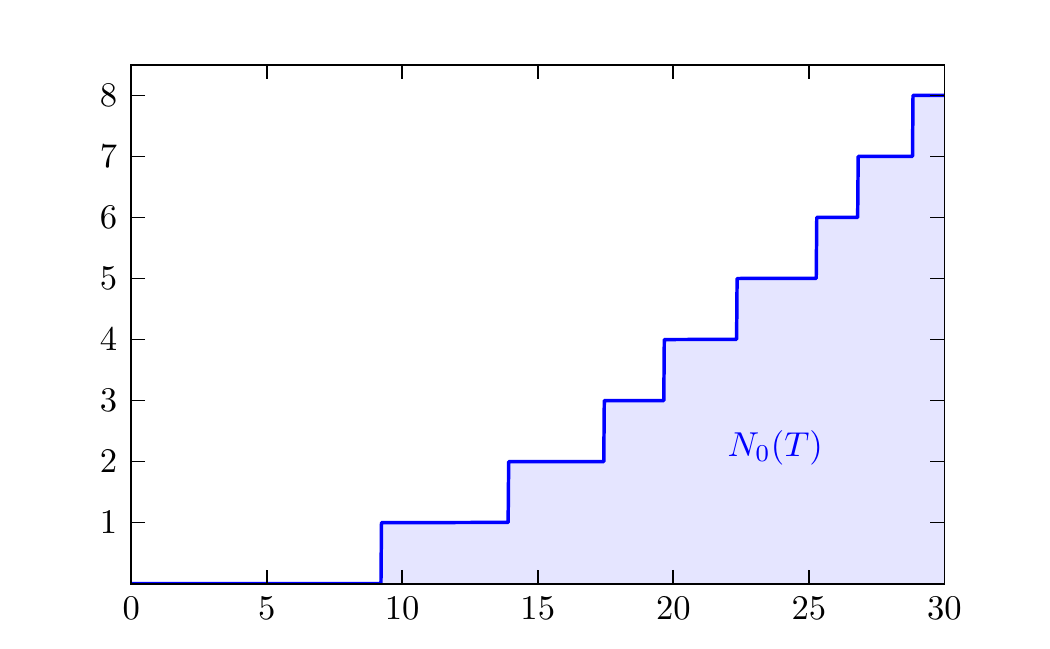}
\caption{Exact counting formula \eqref{countRam} based on the
Ramanujan $\tau$-function.}
\label{fig:ram_counting}
\end{figure}

Here we will consider an example of a modular form of weight $k=12$.
The simplest example is based on the Dedekind $\eta$-function
\beq \label{eta}
\eta (\tau) = q^{1/24} \, \prod_{n=1}^\infty (1-q^n),  \qquad
q=e^{2 \pi i \tau}.
\eeq
Up to a simple factor,  $\eta$ is the inverse of the chiral  partition 
function of the free boson conformal field theory \cite{CFT},  where $\tau$ is
the modular parameter of the torus.   
The modular discriminant 
\beq \label{Delta} 
\Delta (\tau ) = \eta (\tau )^{24}  =  \sum_{n=1}^\infty \, \tau(n) \, q^n
\eeq
is a weight $k=12$ modular form.   
It is closely related to the inverse of the partition function 
of the bosonic string in $26$ dimensions,  where $24$ is the number 
of light-cone degrees of freedom  in $26$ spacetime dimensions \cite{String}.     
The Fourier coefficients  $\tau (n)$  correspond to  the 
Ramanujan $\tau$-function, and the first few are 
\beq
\def\arraystretch{1.2}
\begin{tabular}{@{}c|ccccccccc@{}}
$n$       &$1$ &$2$   &$3$   &$4$     &$5$    &$6$     &$7$      &$8$ \\ 
\hline
$\tau(n)$ &$1$ &$-24$ &$252$ &$-1472$ &$4830$ &$-6048$ &$-16744$ &$84480$ 
\end{tabular}
\eeq
We then define the Dirichlet series
\beq \label{LRam}
L_{\Delta}(s)  =  \sum_{n=1}^\infty  \,  \frac{\tau (n) }{n^s}.
\eeq
Applying \eqref{exactmod} the zeros are $\rho_n = 6 + i t_n$,  
where  $t_n$ satisfies the equation
\beq \label{exactRam}
\vartheta_{12} (t_n) + 
\lim_{\delta \to 0^{+}} \arg L_{\Delta}(6 + \delta + i t_n ) = 
\(n-\tfrac{1}{2}\) \pi.
\eeq
The counting formula \eqref{NTmod} and its asymptotic approximation are
\begin{align} 
\label{countRam}
N_0(T) &= \inv{\pi} \vartheta_{12} (T) + 
\inv{\pi} \arg L_{\Delta} (6 + i T)  \\
&\simeq
\dfrac{T}{\pi} \log \( \dfrac{T}{2 \pi e} \)  + \inv{\pi}\arg L_\Delta(6+iT) +
\dfrac{11}{4}.
\end{align}
A plot of \eqref{countRam} is shown in 
Figure~\ref{fig:ram_counting}, and we can see that it is a perfect staircase
function.
The approximate solution \eqref{yLambertmod} now has  the form
\beq \label{LambertRam} 
\tilde{t}_n = \dfrac{ \(n- \tfrac{13}{4} \) \pi }{W\[(2e)^{-1}
\(n-\tfrac{13}{4}\)\]}
\qquad \mbox{for $n=2,3,\dotsc$}.
\eeq
Note that the above equation is valid for $n > 1$ since $W(x)$ is not defined
for $x < -1/e$.

We follow exactly the same numerical procedure, previously discussed in
Section~\ref{sec:numerical_exact} 
and implemented in Appendix~\ref{sec:mathematica}, to solve
the equation \eqref{exactRam} starting with the approximation
provided by \eqref{LambertRam}. Some of these solutions are shown
in Table~\ref{zerosRam} in Appendix~\ref{sec:tables_ramanujan} 
and are accurate to $50$ decimal places.

\section{Concluding remarks}
\label{sec:conclusion}

In this paper we considered non-trivial zeros of the Riemann $\zeta$-function,
Dirichlet $L$-functions and $L$-functions based on level one modular
forms. The same approach was applied to all these classes of functions, 
showing that there are an infinite number of  zeros on the critical line 
in one-to-one correspondence with the zeros of the cosine 
function \eqref{particular_sol}, 
leading to a transcendental equation satisfied by the ordinate of the $n$-th  
zero.   More specifically, for the
Riemann $\zeta$-function these zeros are solutions to \eqref{exact_eq2}, 
for Dirichlet $L$-functions we have \eqref{exact}, and 
for $L$-functions based on level one modular forms
the ordinates of the zeros must satisfy \eqref{exactmod}. It is important
to stress that these equations were \emph{derived on the critical
line}, without assuming the RH.

The implication of our work for the GRH can be summarized 
as follows.  If the corresponding transcendental equation
has a unique solution for every $n$,  the validity of the GRH would follow.
The explanation is very simple. Suppose the transcendental
equation indeed has a unique solution for every $n$.
Then the zeros obtained from its solutions on 
the \emph{critical line} can be counted,  
since they are enumerated by the integer $n$,  
yielding the counting function $N_0 (T)$.   The number of solutions
saturate the counting formula over the entire 
\emph{critical strip},  namely $N_0 (T) = N(T)$,  
where $N(T)$ counts zeros on the entire critical strip and has been 
known for a long time.   
Thus the equation captures all the non-trivial zeros.

As previously discussed, and explicitly illustrated
in Section~\ref{sec:davenport}, the existence of solutions depends on 
whether the 
$\delta\to0^+$ limit of the argument of the corresponding $L$-function 
is well-defined 
for every ordinate $t$. The validity of this limit was our only assumption
throughout the paper.   
We also argued that if there is indeed a unique solution 
of the transcendental equation for every $n$,  then all non-trivial
zeros are simple.
If there are zeros off of the critical
line, or zeros with multiplicity $m\ge2$ on the critical line, the
equation will fail to capture all the zeros on the critical strip
and then $N_0(T) < N(T)$. Does this means that the GRH is false if
the transcendental equation does not always have a unique solution?
Not necessarily, since all the zeros can still be on the critical line but
not all of them are simple. This suggests that the GRH and the simplicity
of all non-trivial zeros is equivalent
to the statement that the transcendental equation has a unique solution
for every $n$.     We have not proven that there is a unique solution 
to the transcendental equation.    
An attempt to justify more carefully this $\delta\to0^+$
limit is  in our preliminary work \cite{EulerProd},
where we claim that the Euler Product Formula is still valid in the
region $\tfrac{1}{2}< \Re(s) \le 1$ in a statistical manner.

We also showed that it is possible to obtain an explicit formula as  an
approximation for the ordinates of the zeros in terms of the 
Lambert $W$-function; 
equation \eqref{Lambert}
for the $\zeta$-function, \eqref{approx} for Dirichlet $L$-functions
and \eqref{yLambertmod} for $L$-functions based on level one modular forms.
This approximation is very convenient, allowing us to actually compute
accurate zeros without relying on Gram points, nor dealing with violations
of Gram's law.

We have also provided compelling  numerical evidence for the validity of these 
transcendental equations satisfied by the $n$-th zero.
For the $\zeta$-function, the leading order asymptotic approximation
\eqref{FinalTranscendence} proved to be accurate enough to reveal
the interesting features of the Riemann zeros, like the GUE statistics
and the reconstruction of the prime number counting function $\pi(x)$.
It turns out the exact equation \eqref{exact_eq2} is much
more stable and easy to solve numerically,  
it is thus able to provide numerical results as accurate as is desired.
We have  also provided accurate numerical solutions for 
Dirichlet $L$-functions using \eqref{exact} and for the particular
example of the modular $L$-function based on the Ramanujan $\tau$-function,
through \eqref{exactRam}. The numerical approach employed here
constitutes a novel and simple method to compute non-trivial zeros
of $L$-functions.

\section*{Acknowledgments}
We wish to thank Tim Healey and Wladyslaw Narkiewicz for useful 
discussions. We are grateful to the anonymous referee for valuable
suggestions.
GF is supported by CNPq-Brazil.

\appendix

\section{Mathematica implementation}
\label{sec:mathematica}

Here we provide the short Mathematica code used to compute the zeros
from the transcendental equations. We will consider
Dirichlet $L$-functions, since it involves more ingredients, like the
modulus $k$, the order $a$ and the Gauss sum $G(\tau)$. For the Riemann
$\zeta$-function the procedure below is trivially adapted as a special case, as
well as for the Ramanujan $\tau$-function of Section~\ref{sec:ramanujan}.

The function \eqref{RSgen} is implemented as follows:
\begin{lstlisting}
RSTheta[t_, a_, k_] := Im[LogGamma[1/4+a/2+I*t/2]] - t/2*Log[Pi/k]
\end{lstlisting}
For the transcendental equation \eqref{exact} we have
\begin{lstlisting}
ExactEq[n_, t_, s_, a_, k_, j_, G_, n0_] := 
    (RSTheta[t, a, k] + Arg[DirichletL[k, j, 1/2+\delta+I*t]] - 1/2*Arg[G])/Pi + a/4 + 1/2 - n + n0
\end{lstlisting}
Above, $s$ denotes $0<\delta\ll 1$, $a$ is the order \eqref{order}, $k$ is the
modulus, $j$ specifies the Dirichlet character $\chi_{k,j}$ (as discussed
in Section~\ref{sec:dirichlet}), and $G$ is the Gauss
sum \eqref{tau}. Note that we also included $n_0$, discussed
after \eqref{almost_final}, but we always set $n_0=0$ for the cases
analysed in Section~\ref{sec:numerical_dirichlet}. The implementation of the
approximate solution \eqref{approx} is
\begin{lstlisting}
Sgn[n_] := Which[n != 0, Sign[n], n == 0, -1]
A[n_, a_, G_, n0_] := Sgn[n]*(n - n0 + 1/2/Pi*Arg[G]) 
        + (1 - 4*Sgn[n] - 2*a*(Sgn[n]+1))/8
tApprox[n_, a_, G_, k_, n0_] := 
    2*Pi*Sgn[n]*A[n, a, G, n0]/LambertW[k*A[n, a, G, n0]/E]
\end{lstlisting}
One can then obtain the numerical solution of the transcendental 
equation \eqref{exact} as follows:
\begin{lstlisting}
FindZero[n_, s_, a_, k_, j_, G_, n0_, t0_, prec_] := 
    t /. FindRoot[ExactEq[n, t, s, a, k, j, G, n0], {t, t0}, PrecisionGoal->prec/2, AccuracyGoal->prec/2, WorkingPrecision->prec]
\end{lstlisting}
Above, $t_0$ will be given by the approximate solution \eqref{approx}.
The variable $prec$ will be adjusted iteratively.
Now the procedure described in Section~\ref{sec:numerical_exact}
can be implemented as follows:
\begin{lstlisting}
DirichletNZero[n_, order_, digits_, k_, j_, n0_] := (
    chi = DirichletCharacter[k, j, -1];
    a = Which[chi == -1, 1, chi == 1, 0];
    s = 10^(-3);
    prec = 15;
    G = Sum[DirichletCharacter[k, j, l]*Exp[2*Pi*I*l/k], {l, 1, k}];
    t = N[tApprox[n, a, G, k, n0], 20];
    absvalue = 1;
    While[absvalue > order,
        t = FindZero[n, s, a, k, j, G, n0, t, prec];
        Print[NumberForm[t, digits]];
        s = s/1000;
        prec = prec + 20;
        absvalue = Abs[DirichletL[k, j, 1/2 + I*t]];
    ]
    Print[ScientificForm[absvalue, 5]];
)
\end{lstlisting}
Above the variable $order$ controls the accuracy of the solution. 
For instance, if $order=10^{-50}$ it will iterate until 
$|L\(\tfrac{1}{2})+it\)| \sim 10^{-50}$. The
variable $digits$ controls the number of decimal places shown in the output.

Let us compute the zero $n=1$, for the character \eqref{char73},
i.e. $k=7$ and $j=3$. We will verify the solution to $order=10^{-20}$ and
print the results with $digits=22$. Thus executing
\begin{lstlisting}
DirichletNZero[1, 10^(-20), 22, 7, 3, 0]
\end{lstlisting}
the output will be
\begin{lstlisting}
4.35640188194944
4.356401624736541498075
4.356401624736284227537
4.356401624736284227280
4.1664*10^(-25)
\end{lstlisting}
Note how the decimal digits converge in each iteration. 
By decreasing $order$ and
increasing $digits$ it is possible to obtain  highly accurate solutions.
It is exactly in this way that we obtained the tables shown in
Appendix~\ref{sec:tables}. Obviously, depending on the height of
the critical line under consideration, one should adapt the 
parameters $s$ and $prec$ appropriately. In Mathematica we were able to
compute solutions up to  $n\sim10^6$ for Dirichlet $L$-functions, and up to 
$n\sim10^9$ for the Riemann $\zeta$-function without problems.  
We were unable to go much higher only because Mathematica could not 
compute the $\arg L$ term reliably.     
 To solve the transcendental
equations \eqref{exact_eq2} and \eqref{exact} for very high values
on the critical line is still a  challenging numerical problem. Nevertheless,
we believe that it can be done through a more specialized implementation.

\section{Numerical results}
\label{sec:tables}

In this section we present some of the numerical results obtained by
solving the transcendental equations described in this paper. The
numerical procedure is described in Appendix~\ref{sec:mathematica} and
should be adapted to each particular class of functions.

\subsection{Riemann $\zeta$-function}
\label{sec:tables_zeta}

The explicit formula \eqref{Lambert} can estimate very high Riemann zeros,
yielding results accurate up to the decimal point. Some of these results are
shown below:
\begin{table}[H]
\def\arraystretch{1}
\centering
\footnotesize
\begin{tabular}{@{}ll@{}}
\toprule[0.8pt]
$n$ & $\tilde{t}_n$ \\
\midrule[0.4pt]
$10^{22}$ & $1.370919909931995308226636$ \hfill $\times10^{21}$ \\
$10^{50}$ & $5.741532903784313725642221053588442131126693322343461$ 
  \hfill $\times10^{48}$ \\ 
          &  \hfill \footnotesize{\it(continued)}
\end{tabular}
\end{table}
\begin{table}[H]
\centering
\def\arraystretch{1}
\footnotesize
\begin{tabular}{@{}rrr@{}}
            &  \hfill \footnotesize{\it(continued)} \\
$10^{100}$  & 
$2.806903838428940699031954458382564000845480301628460$\\[-.4ex]
  &\phantom{$2.$}$45192360059224930922349073043060335653109252473234$
  \hfill $\times10^{98}$ \\ 
$10^{200}$ &
$1.385792222146789340845466805467159190123402451538707$\\[-.4ex]
  &\phantom{$1.$}$081832868352483938909689796343076797639408172610028$\\[-.4ex]
  &\phantom{$1.$}$651791994879400728026863298840958091288304951600695$\\[-.4ex]
  &\phantom{$1.$}$814960962282888090054696215023267048447330585768$
  \hfill$\times10^{198}$ \\
\bottomrule[0.8pt]
\end{tabular}
\caption{Numerical results predicted by formula \eqref{Lambert}, which 
can easily estimate very high Riemann zeros.
The results are expected to be correct up to the decimal point,  i.e. to the 
number of digits in the integer part.   
The numbers are shown with three digits beyond the integer part.}
\label{highn}
\end{table}

In the following table we have numerical solutions to 
\eqref{FinalTranscendence}, obtained simply by applying a root
finder method around the estimate provided by formula \eqref{Lambert}:
\begin{table}[H]
\def\arraystretch{1}
\centering
\footnotesize
\begin{tabular}{@{}lrr@{}}
\toprule[0.8pt]
$n$ & $\tilde{t}_n$ & $t_n$ \\ 
\midrule[0.4pt]
$1$       &         $14.52$ &         $14.134725142$ \\
$10$      &         $50.23$ &         $49.773832478$ \\
$10^2$    &        $235.99$ &        $236.524229666$ \\
$10^3$    &       $1419.52$ &       $1419.422480946$ \\
$10^4$    &       $9877.63$ &       $9877.782654006$ \\
$10^5$    &      $74920.89$ &      $74920.827498994$ \\
$10^6$    &     $600269.64$ &     $600269.677012445$ \\
$10^7$    &    $4992381.11$ &    $4992381.014003179$ \\
$10^8$    &   $42653549.77$ &   $42653549.760951554$ \\
$10^9$    &  $371870204.05$ &  $371870203.837028053$ \\
$10^{10}$ & $3293531632.26$ & $3293531632.397136704$ \\
\bottomrule[0.8pt]
\end{tabular}
\caption{Numerical solutions to the asymptotic 
equation \eqref{FinalTranscendence}.
All numbers shown are accurate up to the $9$-th decimal place and
agree with \cite{Odlyzko}.}
\label{some_zeros}
\end{table}

Decreasing the error tolerance we can obtain more accurate solutions 
to the asymptotic equation \eqref{FinalTranscendence}, even for the 
lower zeros, as shown below:
\begin{table}[H]
\def\arraystretch{1.2}
\centering
\footnotesize
\begin{tabular}{@{}ll@{}}
\toprule[0.8pt]
$n$ & $t_n$ \\ 
\midrule[0.4pt]
$1$ & $14.13472514173469379045725198356247$ \\
$2$ & $21.02203963877155499262847959389690$ \\
$3$ & $25.01085758014568876321379099256282$ \\
$4$ & $30.42487612585951321031189753058409$ \\
$5$ & $32.93506158773918969066236896407490$ \\
\bottomrule[0.8pt]
\end{tabular}
\caption{Numerical solutions to \eqref{FinalTranscendence}  for  the 
lowest zeros.   Although it  was derived for high $t$,
it provides accurate solutions even for the lower zeros. 
These  numbers are correct up to the decimal place shown \cite{Odlyzko}.}
\label{lower_zeros}
\end{table}

While the previous tables were computed for isolated zeros,
to test Odlyzko-Montgomery pair correlation conjecture \eqref{odlyzko_pair}
we have to compute systematically a wide range of zeros. This is
a strong test of equation \eqref{FinalTranscendence} and 
the approximation \eqref{Lambert}, 
since in principle it could have missed some zeros or presented
some numerical issues. This was definitely not the case.
We computed all the zeros in the range $n=1\dotsc10^5$ and 
also $n=10^9-10^5\dotsc10^9$. The equation \eqref{FinalTranscendence}
and also the approximation \eqref{Lambert} did not miss a single zero.
The last numbers in these ranges are shown below:
\begin{table}[H]
\def\arraystretch{1}
\centering
\footnotesize
\begin{minipage}{.49\textwidth}
\centering
\begin{tabular}{@{}ll@{}}
\toprule[0.8pt]
$n$ & $t_n$ \\ 
\midrule[0.4pt]
$10^5-5$ & $74917.719415828$ \\
$10^5-4$ & $74918.370580227$ \\
$10^5-3$ & $74918.691433454$ \\
$10^5-2$ & $74919.075161121$ \\
$10^5-1$ & $74920.259793259$ \\
$10^5$   & $74920.827498994$ \\
\bottomrule[0.8pt]
\end{tabular}
\end{minipage}
\begin{minipage}{.49\textwidth}
\centering
\begin{tabular}{@{}ll@{}}
\toprule[0.8pt]
$n$ & $t_n$ \\ 
\midrule[0.4pt]
$10^9-5$ & $371870202.244870467$ \\
$10^9-4$ & $371870202.673284457$ \\
$10^9-3$ & $371870203.177729799$ \\
$10^9-2$ & $371870203.274345928$ \\
$10^9-1$ & $371870203.802552324$ \\
$10^9$   & $371870203.837028053$ \\
\bottomrule[0.8pt]
\end{tabular}
\end{minipage}
\caption{Last numerical solutions to \eqref{FinalTranscendence} around
$n=10^5$ and $n=10^9$.}
\label{high_values}
\end{table}

The previous numerical solutions to \eqref{FinalTranscendence} were
obtained with no iteration, i.e. simply by applying the root
finder function once.

The numerical solutions to the exact equation \eqref{exact_eq2}
can yield arbitrarily accurate values. With some very few iterations,
as described in Appendix~\ref{sec:mathematica}, we computed the first 
few zeros:
\begin{table}[H]
\centering
\def\arraystretch{1.2}
\footnotesize
\begin{tabular}{@{}ll@{}}
\toprule[0.8pt]
$n$ & $t_n$ \\ 
\midrule[0.4pt]
$1$ & $14.1347251417346937904572519835624702707842571156992431756855$ \\
$2$ & $21.0220396387715549926284795938969027773343405249027817546295$ \\
$3$ & $25.0108575801456887632137909925628218186595496725579966724965$ \\
$4$ & $30.4248761258595132103118975305840913201815600237154401809621$ \\
$5$ & $32.9350615877391896906623689640749034888127156035170390092800$ \\
\bottomrule[0.8pt]
\end{tabular}
\caption{The first few numerical solutions to \eqref{exact_eq2},
accurate to $60$ digits ($58$ decimals).}
\label{lower_precise}
\end{table}

\subsection{Dirichlet $L$-functions}
\label{sec:tables_dirichlet}

Below we present some numerical solutions to \eqref{exact}, with the
Dirichlet character shown in \eqref{char72}. We used exactly the
procedure described in Appendix~\ref{sec:mathematica}. For $n>0$ we have
the zeros on the upper half of the critical line, while for
$n\le 0$ we have the zeros on the lower half of the critical line.
\begin{table}[H]
\centering
\footnotesize
\def\arraystretch{1}
\begin{tabular}{@{}rrr@{}}
\toprule[0.8pt]
$n$ & $\tilde{t}_n$ & $t_n$ \\
\midrule[0.4pt]
$10$ &  $25.57$ &  $25.68439458577475868571703403827676455384372032540097$ \\
$9$ &   $23.67$ &  $24.15466453997877089700472248737944003578203821931614$ \\
$8$ &   $21.73$ &  $21.65252506979642618329545373529843196334089625358303$ \\
$7$ &   $19.73$ &  $19.65122423323359536954110529158230382437142654926200$ \\
$6$ &   $17.66$ &  $17.16141654370607042290552256158565828745960439000612$ \\
$5$ &   $15.50$ &  $15.74686940763941532761353888536874657958310887967059$ \\
$4$ &   $13.24$ &  $13.85454287448149778875634224346689375234567535103602$ \\
$3$ &   $10.81$ &   $9.97989590209139315060581291354262017420478655402522$ \\
$2$ &    $8.14$ &   $8.41361099147117759845752355454727442365106861800819$ \\
$1$ &    $4.97$ &   $5.19811619946654558608428407430395403442607551643259$ \\
    &            &  \hfill \footnotesize{\it(continued)}
\end{tabular}
\end{table}
\begin{table}[H]
\centering
\def\arraystretch{1}
\footnotesize
\begin{tabular}{@{}rrr@{}}
    &            &  \hfill \footnotesize{\it(continued)} \\
$0$ &   $-3.44$ &  $-2.50937455292911971967838452268365746558148671924805$ \\
$-1$ &  $-7.04$ &  $-7.48493173971596112913314844807905530366284046079242$ \\
$-2$ &  $-9.85$ &  $-9.89354379409772210349418069925221744973779313289503$ \\
$-3$ & $-12.35$ & $-12.25742488648921665489461478678500208978360618268664$ \\
$-4$ & $-14.67$ & $-14.13507775903777080989456447454654848575048882728616$ \\
$-5$ & $-16.86$ & $-17.71409256153115895322699037454043289926793578042465$ \\
$-6$ & $-18.96$ & $-18.88909760017588073794865307957219593848843485334695$ \\
$-7$ & $-20.99$ & $-20.60481911491253262583427068994945289180639925014034$ \\
$-8$ & $-22.95$ & $-22.66635642792466587252079667063882618974425685038326$ \\
$-9$ & $-24.87$ & $-25.28550752850252321309973718800386160807733038068585$ \\
\bottomrule[0.8pt]
\end{tabular}
\caption{Numerical solutions to \eqref{exact} starting
with the approximation \eqref{approx}, for the
character \eqref{char72}. The solutions are
accurate to $50$ decimal places.}
\label{zeros_1}
\end{table}

The next two tables contains numerical solutions to \eqref{exact},
but with the Dirichlet character \eqref{char73}.
\begin{table}[H]
\centering
\def\arraystretch{1}
\footnotesize
\begin{tabular}{@{}rrr@{}}
\toprule[0.8pt]
$n$ & $\tilde{t}_n$ & $t_n$ \\
\midrule[0.4pt]
$10$ &  $25.55$  &  $26.16994490801983565967242517629313321888238615283992$ \\
$9$ &   $23.65$  &  $23.20367246134665537826174805893362248072979160004334$ \\
$8$ &   $21.71$  &  $21.31464724410425595182027902594093075251557654412326$ \\
$7$ &   $19.71$  &  $20.03055898508203028994206564551578139558919887432101$ \\
$6$ &   $17.64$  &  $17.61605319887654241030080166645399190430725521508443$ \\
$5$ &   $15.48$  &  $15.93744820468795955688957399890407546316342953223035$ \\
$4$ &   $13.21$  &  $12.53254782268627400807230480038783642378927939761728$ \\
$3$ &   $10.79$  &  $10.73611998749339311587424153504894305046993275660967$ \\
$2$ &    $8.11$  &   $8.78555471449907536558015746317619235911936921514074$ \\
$1$ &    $4.93$  &   $4.35640162473628422727957479051551913297149929441224$ \\
$0$ &   $-5.45$  &  $-6.20123004275588129466099054628663166500168462793701$ \\
$-1$ &  $-8.53$  &  $-7.92743089809203774838798659746549239024181788857305$ \\
    &            &  \hfill \footnotesize{\it(continued)}
\end{tabular}
\end{table}

\begin{table}[H]
\centering
\def\arraystretch{1}
\footnotesize
\begin{tabular}{@{}rrr@{}}
    &            &  \hfill \footnotesize{\it(continued)} \\
$-2$ & $-11.15$  & $-11.01044486207249042239362741094860371668883190429106$ \\
$-3$ & $-13.55$  & $-13.82986789986136757061236809479729216775842888684529$ \\
$-4$ & $-15.80$  & $-16.01372713415040781987211528577709085306698639304444$ \\
$-5$ & $-17.94$  & $-18.04485754217402476822077016067233558476519398664936$ \\
$-6$ & $-20.00$  & $-19.11388571948958246184820859785760690560580302023623$ \\
$-7$ & $-22.00$  & $-22.75640595577430793123629559665860790727892846161121$ \\
$-8$ & $-23.94$  & $-23.95593843516797851393076448042024914372113079309104$ \\
$-9$ & $-25.83$  & $-25.72310440610835748550521669187512401719774475488087$ \\
\bottomrule[0.8pt]
\end{tabular}
\caption{Numerical solutions of \eqref{exact} starting
with the approximation \eqref{approx}, for the character \eqref{char73}. 
The solutions are
accurate to $50$ decimal places.}
\label{zeros_2}
\end{table}

\begin{table}[H]
\def\arraystretch{1}
\centering
\footnotesize
\begin{tabular}{@{}crr@{}}
\toprule[0.8pt]
$n$ & $\tilde{t}_n$ & $t_n$ \\
\midrule[0.4pt]
$10^3$ & $1037.61$ & 
  $1037.563717069206542965600461276981687171127496013595490$ \\
$10^4$ & $7787.18$ &      
  $7787.337916840954922060149425635486826208937584171726906$ \\
$10^5$ & $61951.04$ &  
 $61950.779420880674657842482173403370835983852937763461400$ \\
$10^6$ & $512684.78$ & 
$512684.856698029779109684519709321053301710419463624401290$ \\
\bottomrule[0.8pt]
\end{tabular}
\caption{Higher zeros for the Dirichlet character \eqref{char73}.  
These solutions to \eqref{exact}  are accurate to $50$ decimal places.}
\label{highzeros}
\end{table}

\subsection{$L$-function based on Ramanujan $\tau$}
\label{sec:tables_ramanujan}

Adapting the numerical procedure of \ref{sec:mathematica} for the
modular $L$-function based on the Ramanujan $\tau$-function, described
in section \ref{sec:ramanujan}, we can obtain the following
numerical solutions, some of which were previously unknown:
\begin{table}[H]
\def\arraystretch{1}
\centering
\footnotesize
\begin{tabular}{@{}rrr@{}}
\toprule[0.8pt]
$n$ & $\tilde{t}_n$ & $t_n$ \\
\midrule[0.4pt]
$1$   &          &   $9.22237939992110252224376719274347813552877062243201$ \\
$2$   &  $12.46$ &  $13.90754986139213440644668132877021949175755235351449$ \\
    &            &  \hfill \footnotesize{\it(continued)}
\end{tabular}
\end{table}
\begin{table}[H]
\centering
\def\arraystretch{1}
\footnotesize
\begin{tabular}{@{}rrr@{}}
    &            &  \hfill \footnotesize{\it(continued)} \\
$3$   &  $16.27$ &  $17.44277697823447331355152513712726271870886652427527$ \\
$4$   &  $19.30$ &  $19.65651314195496100012728175632130280161555091200324$ \\
$5$   &  $21.94$ &  $22.33610363720986727568267445923624619245504695246527$ \\
$6$   &  $24.35$ &  $25.27463654811236535674532419313346311859592673122941$ \\
$7$   &  $26.60$ &  $26.80439115835040303257574923358456474715296800497933$ \\
$8$   &  $28.72$ &  $28.83168262418687544502196191298438972569093668609124$ \\
$9$   &  $30.74$ &  $31.17820949836025906449218889077405585464551198966267$ \\
$10$  &  $32.68$ &  $32.77487538223120744183045567331198999909916163721260$ \\
$100$ & $143.03$ & $143.08355526347845507373979776964664120256210342087127$ \\
$200$ & $235.55$ & $235.74710143999213667703807130733621035921210614210694$ \\
$300$ & $318.61$ & $318.36169446742310747533323741641236307865855919162340$ \\
\bottomrule[0.8pt]
\end{tabular}
\caption{Non-trivial zeros of the modular $L$-function based
on the Ramanujan $\tau$-function, obtained from \eqref{exactRam} starting
with the approximation \eqref{LambertRam}.  
These solutions are accurate to $50$ decimal places.}
\label{zerosRam}
\end{table}



\begin{thebibliography}{99}

\bibitem{Riemann} B. Riemann,
{\it Ueber die Anzahl der Primzahlen unter einer gegebenen Gr\" osse},
Monatsberichte der Berliner Akademie. In Gesammelte Werke, Teubner, Leipzig
(1982). See english translation in \cite{Edwards}, {\it ``On the 
number of primes less than a given magnitude''}

\bibitem{Edwards}    H. M. Edwards, 
{\it Riemann's Zeta Function}, Dover Publications Inc., 1974

\bibitem{Sarnak} P. Sarnak, {\it Problems of the Millennium: The Riemann
hypothesis}, Clay Mathematics Institute (2004)

\bibitem{Bombieri} E. Bombieri, {\it Problems of the Millennium: The
Riemann hypothesis}, Clay Mathematics Institute (2000)

\bibitem{Conrey}  J.  B.  Conrey,  
{\it The Riemann Hypothesis},
Notices of the AMS {\bf 50} (2003) 342  

\bibitem{Levinson} N. Levinson, {\it More than one third of the zeros
of Riemann's zeta-function are on $\sigma=1/2$}, Advances in Math. {\bf 13}
(1974) 383--436

\bibitem{Conrey2} J. B. Conrey, {\it More than two fifths of the zeros
of the Riemann zeta function are on the critical line},
J. reine angew. Math. {\bf 399} (1989) 1--26

\bibitem{Bui} H. Bui, J. B. Conrey, M. Young,
{\it More than $41\%$ of the zeros of the zeta function are on the critical
line}, Acta Arith. 150 (2011) 35--64

\bibitem{Feng}S. Feng, {\it Zeros of the Riemann zeta function on
the critical line}, arXiv:1003.0059 [math.NT] (2010)

\bibitem{Titchmarsh} E. C. Titchmarsh, {\it The Theory of the Riemann
Zeta-Function}, Oxford University Press, 1988

\bibitem{Schumayer} D. Schumayer, D. A. W. Hutchinson,
{\it Physics of the Riemann hypothesis},
Rev. Mod. Phys. {\bf 83} (2011) 307--330

\bibitem{Julia}B. Julia, {\it Statistical theory of numbers},
Number Theory and Physics, eds. J. M. Luck, P. Moussa, and
M. Waldshmidt, Springer Proceedings in Physics, Vol. {\bf 47},
Springer-Verlag, 1990, pp. 276-293

\bibitem{Spector} D. Spector, {\it Supersymmetry and the M\" obius
inversion function}, Comm. Math. Phys. {\bf 127} (1990) 239--252

\bibitem{Bakas}I. Bakas, M. J. Bowick,
{\it Curiosities of arithmetic gases}, J. Math. Phys. {\bf 32} (1991) 1881

\bibitem{Spector2}D. Spector, {\it Duality, partial supersymmetry,
and arithmetic number theory}, J. Math. Phys. {\bf 39} (1998) 1919--1927

\bibitem{Cacciatori} S. L. Cacciatori, M. Cardella,
{\it Equidistribution rates, closed string amplitudes, and the Riemann
hypothesis}, JHEP {\bf 12} (2010)

\bibitem{He} Y-H. He, V. Jejjalla, D. Minic,
{\it On the physics of the Riemann zeros},
J. Phys: Conf. Ser. {\bf 462} (2013) 012036

\bibitem{Lapidus}M. L. Lapidus,
{\it In Search of the Riemann Zeros: Strings, Fractal 
Membranes and Noncommutative Spacetimes}, Amer. Math. Soc. (2008)

\bibitem{AL}  A. LeClair, 
{\it Interacting Bose and Fermi gases in low dimensions and the Riemann
hypothesis}, 
Int. J. Mod. Phys. A {\bf 23} (2008) 1371

\bibitem{Montgomery}  H. Montgomery,  
{\it The pair correlation of zeros of the zeta function}, 
Analytic number theory, Proc. Sympos. Pure Math. XXIV, Providence, R.I.: 
AMS, pp. 181D193,  1973

\bibitem{Berry1}   M. V.  Berry,
{\it Riemann's Zeta Function:  A Model for Quantum Chaos?},
Quantum Chaos and Statistical Nuclear Physics, Eds. T. H. Seligman and
H. Nishioka, Lecture Notes in Physics, {\bf 263} Springer Verlag, New
York, 1986

\bibitem{BerryKeating}  M. V. Berry, J. P. Keating,    
{\it The Riemann zeros and eigenvalue asymptotics}, 
SIAM Review {\bf 41} (1999)  236 

\bibitem{BerryKeating2} M. V. Berry, J. P. Keating,  
{\it $H=xp$ and the Riemann zeros}, 
in {\it Supersymmetry and Trace Formulae:  Chaos and Disorder}, Kluwer 1999  

\bibitem{Sierra}   G. Sierra,   
{\it The Riemann zeros and the cyclic renormalization group}, 
J. Stat. Mech. 0512:P12006 (2005)  

\bibitem{Sierra2}    G.  Sierra, 
{\it A Quantum Mechanical model of the Riemann Zeros},
New J. Phys. {\bf 10}  (2008) 033016

\bibitem{Bhaduri}   R. K. Baduri,  A. Khare, J. Law, 
{\it The phase of the Riemann zeta function and the inverted 
harmonic oscillator}, 
Phys. Rev. {\bf E52} (1995)  486

\bibitem{Connes} A. Connes, 
{\it Trace formula in noncommutative geometry and the zeros of the Riemann 
zeta function}, Sel. Math. New Ser. {\bf 5} (1999) 29--106

\bibitem{Connes2} J-B. Bost, A. Connes,
{\it Hecke algebras, type III factors and phase transitions with spontaneous
symmetry breaking in number theory}, Sel. Math. New Ser. 
{\bf 3} (1995) 411--457

\bibitem{Svaiter} G. Menezes, B. F. Svaiter, N. F. Svaiter,
{\it Riemann zeta zeros and prime number spectra in quantum field theory}, 
Int. J. Mod. Phys {\bf A28} (2013) 1350128

\bibitem{Andrade} J. C. Andrade,
{\it Hilbert-P\' olya conjecture, zeta functions and bosonic quantum field
theories}, Int. J. Mod. Phys {\bf A28} (2013) 1350072

\bibitem{Apostol} T. M. Apostol, {\it Introduction to Analytic Number theory},
Springer-Verlag New York, 1976

\bibitem{Selberg1} A. Selberg, {\it Contributions to the theory of 
Dirichlet's $L$-functions}, Skr. Norske Vid. Akad. Oslo. I. {\bf 1946} (1946)
2--62

\bibitem{Fujii} A. Fujii, {\it On the zeros of Dirichlet $L$-functions. I},
Transactions of the American Math. Soc. {\bf 196} (1974) 225--235

\bibitem{Iwaniec} H. Iwaniec, W. Luo, P. Sarnak, {\it Low lying zeros
of families of $L$-functions}, 
Publications Mathématiques de l'Institut des Hautes Études Scientifiques
{\bf 91} 1 (2000) 55--131

\bibitem{Conrey3} J. B. Conrey, {\it $L$-functions and random matrices},
Mathematics unlimited—-2001 and beyond, 331–-352, Springer 2001

\bibitem{Hughes} C. P. Hughes, Z. Rudnick,
{\it Linear statistics of low-lying zeros of $L$-functions},
Quart. J. Math. {\bf 54} (2003) 309--333

\bibitem{Bogomolny} E. B. Bogomolny and J. P. Keating,
{\it Two-point correlation function for Dirichlet $L$-functions},
J. Phys. A {\bf 46} (2013)

\bibitem{Conrey4} J. B. Conrey, H. Iwaniec, K. Soundararajan,
{\it Critical zeros of Dirichlet $L$-functions},
Journal f\" ur die reine und angewandte Mathematik (Crelles Journal) 
{\bf 681} (2013) 175–-198

\bibitem{Bombieri2} E. Bombieri,
{\it The classical theory of Zeta and $L$-functions},
Milan J. Math. {\bf 78} (2010) 11--59

\bibitem{Iwaniec2} H. Iwaniec, E. Kowalski,
{\it Analytic Number Theory},
American Mathematical Society, 2004

\bibitem{Sarnak2} H. Iwaniec, P. Sarnak,
{\it Perspectives on the analytic theory of $L$-functions},
GAFA, Geom. funct. anal. Special Volume (2000) 705--741

\bibitem{deLune} J. van de Lune, H. J. J. te Riele, D. T. Winter,
{\it On the zeros of the Riemann zeta function in the critical strip. IV},
Math. Comp. {\bf 46} (1986) 667--681

\bibitem{Gourdon}    X.  Gourdon, 
{\it The $10^{13}$ first zeros of the Riemann Zeta function, and zeros 
computation at very large height} (2004)

\bibitem{RHLeclair} A. LeClair,
{\it An electrostatic depiction of the validity of the Riemann 
hypothesis and a formula for the $N$-th zero at large $N$}, 
Int. J. Mod. Phys. A {\bf 28} (2013) 1350151

\bibitem{Goldston}D. A. Goldston, S. M. Gonek, {\it A note on $S(t)$ and
the zeros of the Riemann zeta-function},
Bull. London Math. Soc. {\bf 39} (2007) 482--486

\bibitem{Borwein}J. M. Borwein, D. M. Bradley, R. E. Crandall,
{\it Computational strategies for the Riemann zeta function},
J. Comp. and App. Math. {\bf 121} (200) 247--296

\bibitem{Corless}R. M. Corless, G. H. Gonnet, D. E. G. Hare, D. J. Jeffrey, 
D. E. Knuth, 
{\it On the Lambert $W$ function},
Adv. Comp. Math. {\bf 5} (1996) 329--359

\bibitem{Odlyzko} A. M. Odlyzko, 
{\it Tables of zeros of the Riemann zeta function} \\
www.dtc.umn.edu/{~}odlyzko/zeta-tables/

\bibitem{OdlyzkoSchonhage} A.  M. Odlyzko, A. Sch\" onhage, 
{\it Fast algorithms for multiple evaluations of the Riemann zeta function}, 
Trans. Amer. Math. Soc. {\bf 309} (2) (1988): 797D809

\bibitem{Odlyzko2}   A.  M.  Odlyzko,     
{\it  On the Distribution of spacings  between zeros of the zeta Function}, 
Math. of Computation {\bf 177}  (1987)  273

\bibitem{BombieriGosh} E. Bombieri, A. Ghosh,
{\it Around the Davenport-Heilbronn function},
Russ. Math. Surv. {\bf 66} (2011) 221

\bibitem{Lectures} G. Fran\c ca, A. LeClair,
{\it A theory for the zeros of Riemann zeta and other $L$-functions},
Lectures delivered at ``Riemann Master School on Zeta Functions'' at 
Leibniz Universitat Hannover/Germany, June 10-14, 2014,
arXiv:1407.4358 [math.NT]

\bibitem{Apostolmodular} T. M. Apostol,
{\it Modular Functions and Dirichlet Series in Number Theory},
Springer-Verlag, 1990

\bibitem{CFT} P. Francesco, P. Mathieu, D. Senechal,
{\it Conformal Field Theory}, Springer, 1999

\bibitem{String} M. B. Green, J. H. Schwarz, E. Witten,
{\it Superstring Theory: Volume 1}, Cambridge University Press, 1988


\bibitem{Dyson}  F.  Dyson,    
{\it Correlations between eigenvalues of a random matrix}, 
Comm.  Math.  Phys.  {\bf 19}  (1970) 235

\bibitem{Boros} G. Boros, V. H. Moll, {\it Irresistible Integrals},
Cambridge University Press, 2004, pp. 204



\bibitem{EulerProd} G. Fran\c ca, A. LeClair,
{\it On the validity of the Euler product inside the critical strip},
arXiv:1410.3520 [math.NT] (2014)

\end{thebibliography}
\end{document}